\documentclass[notitlepage]{amsart}

\usepackage{amsfonts}
\usepackage[latin1]{inputenc}
\usepackage[all]{xy} 
\usepackage{latexsym,amssymb,amscd}
\usepackage{amsmath}

\usepackage{graphicx}
\usepackage{fancyhdr}
\usepackage{color}

\usepackage{tikz}
\usetikzlibrary{arrows}
 \usetikzlibrary{snakes}
 \usetikzlibrary{shapes}

\newtheorem{pro}{Proposition}[section]
\newtheorem{teo}[pro]{Theorem}
\newtheorem{defi}[pro]{Definition}
\newtheorem{lem}[pro]{Lemma}
\newtheorem{cor}[pro]{Corollary}
\newtheorem{rk}[pro]{Remark}

\newcommand{\ltt}{\mathcal}

\newcommand{\T}{\mathcal{T}}
\newcommand{\I}{\mathcal{I}}
\newcommand{\Pro}{\mathcal{P}}

\newcommand{\C}{\mathcal{C}}
\newcommand{\Qb}{\mathbf{Q}}
\newcommand{\Yb}{\mathbf{Y}}
\newcommand{\Kb}{\mathbf{K}}
\newcommand{\F}{\mathfrak{F}}
\newcommand{\A}{\mathcal{A}}

\newcommand{\E}{\mathcal{E}}
\newcommand{\Ex}{\mathrm{Ex}}
\newcommand{\X}{\mathcal{X}}
\newcommand{\Y}{\mathcal{Y}}
\newcommand{\Z}{\mathbb{Z}}
\newcommand{\N}{\mathbb{N}}
\newcommand{\D}{\mathrm{D}}
\newcommand{\modu}{\mathrm{mod}}
\newcommand{\Modu}{\mathrm{Mod}}
\newcommand{\proj}{\mathrm{proj}}
\newcommand{\Ext}{\mathrm{Ext}}
\newcommand{\Hom}{\mathrm{Hom}}
\newcommand{\End}{\mathrm{End}}
\newcommand{\Ker}{\mathrm{Ker}}
\newcommand{\Coker}{\mathrm{Coker}}
\newcommand{\Ima}{\mathrm{Im}}

\newcommand{\Tr}{\mathrm{Tr}}
\newcommand{\rad}{\mathrm{rad}}
\newcommand{\Supp}{\mathrm{Supp}}

\newcommand{\add}{\mathrm{add}}
\newcommand{\mini}{\mathrm{min}}
\newcommand{\maxi}{\mathrm{max}}
\newcommand{\fl}{\longrightarrow}

\newenvironment{dem}{\noindent\bf Proof. \rm }{$\ \Box$}
\usepackage{latexsym,amssymb,amscd}
\usepackage{amsmath}

\begin{document}
\title[Homological Systems in triangulated categories ]{Homological systems in triangulated categories}
\author{O. Mendoza, V. Santiago.}
\date{}
\begin{abstract} We introduce the notion of homological systems $\Theta$ for triangulated categories. Homological systems generalize, 
on one hand, the notion of stratifying systems in module categories, and on the other hand, the notion of exceptional sequences in 
triangulated categories. We prove that, attached to the homological system $\Theta,$ there are two standardly stratified algebras 
$A$ and $B,$ which are derived equivalent. Furthermore, it is proved that the category $\F(\Theta),$ of the $\Theta$-filtered 
objects in a triangulated category $\T,$ admits in a very natural way an structure of an exact category, and then 
there are exact equivalences between the exact category $\F(\Theta)$ and the exact categories of the $\Delta$-good modules associated 
to the standardly stratified algebras $A$ and $B.$ Some of the obtained results can be seen also under the light of the cotorsion 
pairs in the sense of Iyama-Nakaoka-Yoshino (see \ref{CotorsionThetaSys1} and \ref{CotorsionThetaSys2}). We recall that cotorsion 
pairs are studied extensively in relation with cluster tilting categories, $t$-structures and co-$t$-structures. 
\end{abstract}
\maketitle
\section{Introduction.}

In \cite{ADL, CPS,Dlab} were introduced the notion of
standardly stratified algebras, generalizing the class of quasi-hereditary algebras. The 
standardly stratified algebras have shown to
be homologically interesting because of their relationship with
tilting theory and relative homological algebra. In order to give a categorification of the standard modules and the 
characteristic tilting module associated with an standardly stratified algebra, Erdmann and S\'aenz
developed the notion of Ext-injective stratifying system \cite{ES}.
In that paper, they generalized the standard modules and the
characteristic tilting module, obtained by Ringel in \cite{R}.\\
Afterwards, Marcos, Mendoza and S\'aenz introduced the notions
of stratifying system and Ext-projective stratifying system, and furthermore, they
proved that all this notions are equivalent to the one given by
Erdmann and S\'aenz \cite{MMS1, MMS2}. In \cite{MMS2},  they were able to
prove that, for a given a stratifying system $(\Theta,\leq)$ in $\mathrm{mod}(A),$ there
exists a module $Q$ such that $B:=\mathrm{End}(Q)^{op}$ is a standardly
stratified algebra, and moreover, there exists an exact equivalence between the
$\Theta$-filtered modules in $\mathrm{mod}(A)$ and the $\Delta$-good modules in $\mathrm{mod}(B).$ 
We remark that the considered order $\leq,$ attached to a stratifying system, is a finite linear one. On the other hand, 
Mendoza, S\'aenz and Xi developed the theory of stratifying systems for the more general case 
of a finite pre-ordered set \cite{MSX}.
\

Triangulated categories have its origin in
algebraic geometry and algebraic topology. This kind of categories
have become relevant in many different areas of mathematics.
Although the axioms of a triangulated category seems to be hard at
first sight, it turns out that many categories are endowed with the
structure of a triangulated category. 
\

In this paper, we develop the concept of homological systems in the
setting of artin triangulated $R$-categories. Throughout this
notes, $\T$ will denote an arbitrary triangulated category and
$[1]:\T \to \T$ its suspension functor. We introduce the notion of a $\Theta$-system (see \ref{ss}) in a
triangulated category $\T$, which is the corresponding generalization of stratifying systems 
in the category of modules over an
algebra $A$. We also state the concept of a $\Theta$-projective
system, and we show that a $\Theta$-system determines a unique
$\Theta$-projective system (see  \ref{sisexproy}). One of the main
results of this paper, Theorem \ref{R3}, says that 
for a given $\Theta$-system in an artin  triangulated
$R$-category, there exist two standardly stratified algebras $A$ and $B;$ and moreover, we have 
triangulated equivalences $\D^b(\F(\Theta))\simeq \D^b(A)$ and $\D^b(\F(\Theta))\simeq \D^b(B).$ 
Furthermore, it is proved the
Theorem \ref{TeoRingel}, which is a generalization to the setting of
triangulated categories of a well-known result obtained by Ringel 
\cite[Theorem 1]{R}. The Theorem \ref{TeoRingel} states that, for a 
given a family of objects $\Theta=\{\Theta(i)\}_{i}^{n}$ belonging to  an
artin triangulated $R$-category $\ltt{T}$ and satisfying that
$\mathrm{Hom}_{\ltt{T}}(\Theta(j),\Theta(i)[1])=0$ if $j\geq i,$ 
the subcategory $\F(\Theta)$ of the $\Theta$-filtered objects in $\ltt{T}$ is functorially finite.
\

The notion of cotorsion pair in a triangulated category was
recently introduced by Iyama-Yoshino \cite{IY} and  Nakaoka
\cite{N}. This notion seems to be important since it unifies the
notions of: (a) $t$-structures \cite{BBD}, (b) co-$t$-structures 
\cite{P} and (c) cluster tilting subcategories \cite{Keller1}. By using the theory of 
homological systems, it is
constructed two canonical cotorsion pairs (see
\ref{CotorsionThetaSys1}) and it is also determined the core of those
cotorsion pairs.
\

As a beautiful application (see Theorem \ref{TeoSES}), of the theory of homological systems, we
showed that if  $\E=(\E_1,\E_2,\cdots,\E_t)$ is a strongly
exceptional sequence in the triangulated category $\T,$ then there
exists an equivalence $\D^b(\F(\E))\simeq \D^b(A)$ as triangulated categories for some quasi-hereditary algebra
$A$. Observe that this result generalize \cite[Theorem 6.2]{Bo}.\\
A similar result (see Theorem \ref{TeoES}) holds for an exceptional sequence $\E$ in the bounded 
derived category $\D^b(\mathcal{H}),$ where $\mathcal{H}$ is a hereditary abelian 
$k$-category.\\
We recall that the notion of exceptional sequences has its
origin from the study of vector bundles on projective spaces (see,
for instance, \cite{Bo,GL}) and that strongly exceptional sequences
appear very often in algebraic geometry and provides a
non-commutative model for the study of
algebraic varieties (see \cite{Bo}).
\

The paper is organized as follows: In section 2, we give some
basic notions and properties of triangulated categories, which will
be used in the rest of the work.\\
In section 3, it is established the concept of artin triangulated
$R$-category and we give some technical results that will be
useful when proving that a $\Theta$-system determines a unique
$\Theta$-projective system.\\
In section 4, we study the subcategory $\F(\Theta)$ of $\Theta$-filtered objects in a
triangulated category $\T$ and it is proved the Theorem \ref{TeoRingel},
which is a generalization of \cite[Theorem 1]{R}.\\
In section 5, we focus our attention into the $\Theta$-systems. It is shown that 
a $\Theta$-system determines a unique
$\Theta$-projective system. We also prove that, for a given 
$\Theta$-projective system, the filtration multiplicity
$[M:\Theta(i)]_{\xi}$ does not depend on the given filtration $\xi.$\\
In section 6, we show that, for a given $\Theta$-projective system
$(\Theta,\Qb,\leq)$ in $\T,$ there exists an equivalence between
$\F(\Theta)$ and the subcategory of the $\Delta$-good modules in $\mathrm{mod}(B),$ for some
standardly stratified algebra $B.$\\
In Section 7, we show that the triangulation in $\T$ induces in a natural way an exact structure in $\F(\Theta),$ and prove 
Theorem \ref{R3}, which is one of the main results of the paper.\\
Finally, in section 8, we provide some examples of homological systems.

\section{Preliminaries}
\

In this paper, $\T$ will be a triangulated category and $[1]:\T\rightarrow\T$ its
 suspension (shift) functor. Moreover, when we say that $\mathcal{C}$ is a subcategory of $\T,$ it always means that $\mathcal{C}$ is  a full subcategory which is additive and closed under isomorphisms. For a class $\X$ of objects of $\T,$ we denote by $\add\,(\X)$ the smallest subcategory of $\T$ containing $\X,$ closed under finite direct sums and direct summands.
\

For some classes $\X$ and $\Y$ of objects in $\T,$ we write
${}^\perp\X:=\{Z\in\T\,:\,\Hom_\T(Z,-)|_{\X}=0\}$ and $\X^\perp:=\{Z\in\T\,:\,\Hom_\T(-,Z)|_{\X}=0\}.$
\noindent We also recall that $\X*\Y$ denotes the class of objects $Z\in\T$ for which there exists a distinguished triangle $X\rightarrow Z\rightarrow Y\rightarrow X[1]$ in $\T$ with $X\in\X$ and $Y\in\Y.$ Furthermore, it is said that  $\X$ is {\bf closed under extensions} if $ \X*\X\subseteq \X.$
\vspace{.4cm}

We will make use of the following constructions in triangulated categories: the base and co-base change. These constructions remind us the pull-back and the push-out, respectively, of short 
exact sequences in abelian categories.

\begin{pro}\label{basecobase} \cite[2.1]{Beli} For any triangulated category $\ltt{T},$ each one 
of the following conditions is equivalent to the octahedral axiom.
\begin{enumerate}
\item [(a)] BASE CHANGE. For any distinguished triangle
$A\stackrel{f}{\to} B\stackrel{g}{\to} C\stackrel{h}{\to}
A[1]$ and any morphism $\epsilon: E\to C$ in $\T,$ there exists a
commutative diagram in $\T$
$$\xymatrix{&  M\ar@{=}[r]\ar[d]^{\alpha} & M\ar[d]^{\delta}\\
A\ar[r]^{f'}\ar@{=}[d] & G\ar[r]^{g'}\ar[d]^{\beta} &
E\ar[r]^{h'}\ar[d]^{\epsilon} & A[1]\ar@{=}[d]\\
A\ar[r]^{f} & B\ar[r]^{g}\ar[d]^{\gamma} &
C\ar[r]^{h}\ar[d]^{\zeta} & A[1]\\
& M[1]\ar@{=}[r] & M[1],}$$ where all the rows and
columns, in the preceding diagram, are distinguished triangles.
\item [(b)] CO-BASE CHANGE.  For any distinguished triangle
$A\stackrel{f}{\to} B\stackrel{g}{\to} C\stackrel{h}{\to}
A[1]$ and any morphism $\alpha: A\to D$ in $\T,$ there exists a
commutative diagram in $\T$
$$\xymatrix{&  N\ar@{=}[r]\ar[d]^{\zeta} & M\ar[d]^{\delta}\\
C[-1]\ar[r]^{-h[-1]}\ar@{=}[d] &
A\ar[r]^{f}\ar[d]^{\alpha} &
B\ar[r]^{g}\ar[d]^{\beta} & C\ar@{=}[d]\\
C[-1]\ar[r]^{-h'[-1]} &
D\ar[r]^{f'}\ar[d]^{\eta} &
F\ar[r]^{g'}\ar[d]^{\vartheta} & C\\
& N[1]\ar@{=}[r] & N[1],}$$ where all the rows and
columns, in the preceding diagram, are distinguished triangles.
\end{enumerate}
\end{pro}

\begin{lem}\label{serpiente} Consider the following commutative diagram in a triangulated category $\ltt{T}$
$$\xymatrix{A\ar[r]^{\alpha}\ar[d]^{\beta} & B\ar[r]\ar[d]^{\beta'} & C\ar[r]  & A[1]\\
A'\ar[r]^{\alpha'} & B'\ar[r] & C'\ar[r] & A'[1],}$$ 
where the rows are distinguished triangles. Then, the preceding diagram can be completed to the following one
$$\xymatrix{A''[-1]\ar[r]\ar[d] &
B''[-1]\ar[r]\ar[d] & C''[-1]\ar[r]\ar[d] &
A''\ar[d]\\
A\ar[r]^{\alpha}\ar[d]^{\beta} & B\ar[r]\ar[d]^{\beta'} &
C\ar[r]\ar[d]^{\Phi} & A[1]\ar[d]\\
A'\ar[r]^{\alpha'}\ar[d] & B'\ar[r]\ar[d] & C'\ar[r]\ar[d]\ar@{}[dr]|{IX} & A'[1]\ar[d]\\
A''\ar[r] & B''\ar[r]  & C''\ar[r] & A''[1],}$$ where the rows
and columns, in the above diagram, are distinguished triangles and all the squares
commute, except by the one marked with IX, which anti-commutes.
\end{lem}
\begin{dem} By completing $\beta$ and $\beta'$ to distinguished triangles, we have the following
commutative diagram in $\T$
$$\xymatrix{A\ar[r]^{\beta}\ar[d]^{\alpha} & A'\ar[r]^{\gamma}\ar[d]^{\alpha'} & A''\ar[r]^{\delta}  & A[1]\\
B\ar[r]^{\beta'} & B'\ar[r]^{\gamma'} & B''\ar[r]^{\delta'} &
B[1].}$$ Then, there is a morphism $h:A'' \to  B''$ in $\T,$ such that the triple $(\alpha,\alpha',h)$ is a morphism of triangles. Hence $h[-1]:A''[-1] \to B''[-1]$ makes commutative the following square
$$\xymatrix{A''[-1]\ar[r]^{h[-1]}\ar[d]_{-\delta[-1]}
& B''[-1]\ar[d]^{-\delta'[-1]}\\
A\ar[r]^{\alpha} & B.}$$ Then, by a Verdier's result (see Exercise 10.2.6, page
378, in \cite{Wei}), we get the lemma.
\end{dem}
\vspace{.4cm}

The following result remind us the so-called Snake's Lemma.

\begin{pro}\label{SnakeLema} Consider the following commutative diagram in a triangulated category $\ltt{T}$
$$\xymatrix{A\ar[r]^{\alpha}\ar[d]^{\beta}\ar@{}[dr] & B\ar[r]\ar[d]^{\beta'} & C\ar[r]\ar[d]^{\beta''}  & A[1]\ar[d]^{\beta[1]}\\
A'\ar[r]^{\alpha'} & B'\ar[r] & C'\ar[r] & A'[1],}$$ where the rows are distinguished triangles. 
If $\Hom_{\ltt{T}}(A,C'[-1])=0$ then the preceding diagram can
be completed to the following one
$$\xymatrix{A''[-1]\ar[r]\ar[d] &
B''[-1]\ar[r]\ar[d] & C''[-1]\ar[r]\ar[d] &
A''\ar[d]\\
A\ar[r]^{\alpha}\ar[d]^{\beta} & B\ar[r]\ar[d]^{\beta'} &
C\ar[r]\ar[d]^{\beta''} & A[1]\ar[d]^{\beta[1]}\\
A'\ar[r]^{\alpha'}\ar[d] & B'\ar[r]\ar[d] & C'\ar[r]\ar[d]\ar@{}[dr]|{IX} & A'[1]\ar[d]\\
A''\ar[r] & B''\ar[r]  & C''\ar[r] & A''[1],}$$ where the rows
and columns, in the above diagram, are distinguished triangles and all the squares
commute, except the one marked with IX, which anti-commutes.
\end{pro}
\begin{dem} Assume that $\Hom_{\ltt{T}}(A,C'[-1])=0.$ By \ref{serpiente}, the square given (in the first diagram) by the morphism $\alpha$ and $\beta$ can be completed to a diagram as
in \ref{serpiente}. We only need to prove that $\Phi=\beta'',$ but this fact follows from 
\cite[Corollary 5]{Gel} page 243, since $\Hom_{\ltt{T}}(A,C'[-1])=0.$
\end{dem}

\begin{defi} Let $\ltt{T}$ be a triangulated category and $\ltt{A}$ be an abelian category. Consider subcategories $\ltt{X}\subseteq \T$ 
and $\ltt{W}\subseteq\ltt{A},$ which are both closed under extensions.  It is said that:
\begin{enumerate}
\item [(a)] a distinguished triangle $\eta:\;A\to B\to  C\to A[1]$ belongs to $\ltt{X}$, that is
$\eta\in \ltt{X}$, if the objects $A,$ $B$ and $C$ belong to $\ltt{X};$

\item [(b)] an additive functor $F:\ltt{X} \to \ltt{W}$ is
$\textbf{exact}$, if for every distinguished triangle $\eta:\;A\to B\to C\to A[1]$ in $\ltt{X}$,
we have that the sequence $F(\eta):\;0\to F(A)\to F(B)\to  F(C)\to  0$ is exact and belongs to $\ltt{W}$.
\end{enumerate}
\end{defi}

We also recall the following well-known definition (see, for example, \cite{B} and \cite{BR}).

\begin{defi} Let $\X$ and $\Y$ be classes of objects in a triangulated category $\T.$ A morphism $f:X\to C$ in $\T$ is said to be an 
{\bf $\X$-precover} of $C$ if $X\in\X$ and $\Hom_\T(X',f):\Hom_\T(X',X)\to\Hom_\T(X',C)$ is surjective $\forall X'\in\X.$ If any $C\in\Y$ admits 
an $\X$-precover, then $\X$ is called a {\bf precovering} class in $\Y.$ By dualizing the definition above, we get the notion of an 
{\bf $\X$-preenvelope} of $C$ and a {\bf preenveloping} class in $\Y.$ Finally, it is said that $\X$ is {\bf{functorially finite}} in $\T$ if $\X$ 
is both precovering and preenveloping in $\T.$
\end{defi}
\vspace{.4cm}

In what follows, we recall some notions and elementary well-known facts about standardly stratified algebras. Let $\Lambda$ be an 
artin $R$-algebra. We denote by $\modu\,(\Lambda)$ to the category of all finitely 
generated left $\Lambda$-modules, and by $\proj\,(\Lambda)$ to the full subcategory of 
$\modu\,(\Lambda)$ whose objects are the projective $\Lambda$-modules. For $M,N\in\modu\,(\Lambda)$, the {\bf trace} $\Tr_M\,(N)$ of $M$ in $N$, 
is the $\Lambda$-submodule of $N$ generated by the images of all
morphisms from $M$ to $N$. For a given natural number $t,$ we set $[1,t]=\{1,2,\cdots, t\}.$
\

We next recall the definition (see \cite{ADL, Dlab, DR, R}) of the class of standard $\Lambda$-modules. Let $n$ be the rank of the Grothendieck group 
$K_0\,(\Lambda)$. We fix a linear order $\leq$ on the set $[1,n]$ and a representative set ${}_\Lambda P=\{{}_{\Lambda}P(i)\;:\;i\in[1,n]\}$ 
containing one module of each iso-class of indecomposable projective $\Lambda$-modules. The set of {\bf standard $\Lambda$-modules} is 
${}_\Lambda\Delta=\{{}_\Lambda\Delta(i):i\in[1,n]\},$ where
${}_\Lambda\Delta(i)={}_{\Lambda}P(i)/\Tr_{\oplus_{j>i}\,{}_{\Lambda}P(j)}\,({}_{\Lambda}P(i))$. Then,
${}_{\Lambda}\Delta(i)$ is the largest factor module of
${}_{\Lambda}P(i)$ with composition factors only amongst
${}_{\Lambda}S(j)$ for $j\leq i.$ 
\

Let $\F({}_{\Lambda}{\Delta})$ be the subcategory of $\modu\,(\Lambda)$ consisting of the
$\Lambda$-modules having a ${}_{\Lambda}{\Delta}$-filtration, that is,
a sequence of submodules $0=M_0\subseteq M_1\subseteq\cdots\subseteq M_s=M$ with factors $M_{i+1}/M_i$
isomorphic to a module in ${}_{\Lambda}{\Delta}$ for all $i$. The pair $(\Lambda,\leq)$ is said to be 
a {\bf standardly stratified algebra}, with respect to the linear
order $\leq$ on the set $[1,n]$, if
$\proj\,(\Lambda)\subseteq\F({}_{\Lambda}{\Delta})$ (see \cite{ADL,CPS, Dlab}). A {\bf quasi-hereditary algebra} is a standardly stratified algebra $(\Lambda,\leq)$ such that 
$\End({}_{\Lambda}{\Delta}(i))$ is a division ring, for each $i\in[1,n].$

\section{triangulated $R$-categories}

Let $R$ be a commutative ring. We recall that an {\bf $R$-category} is a category $\ltt{C}$ satisfying 
the following two conditions: (a) for each pair $X,Y$ of objects in $\C,$ the set of morphisms $\Hom_\C(X,Y)$ is an $R$-module, and (b) the composition of morphisms in $\ltt{C}$ is $R$-bilinear. An $R$-category $\C$ is called {\bf Hom-finite} if $\Hom_\C(X,Y)$ is a finitely generated $R$-module, for each $X,Y\in \C.$  
\

A functor $F:\C\to\ltt{D},$ between $R$-categories, is said to be an {\bf $R$-functor} if $F:\Hom_\C(X,Y)\to \Hom_{\ltt{D}}(F(X), F(Y))$ is a morphism of $R$-modules for each 
pair $X,Y$ of objects in $\C.$

\begin{defi} A {\bf triangulated $R$-category} is an $R$-category $\T$ which is a triangulated one, and 
such that its shift functor $[1]:\T\to\T$ is an $R$-functor.
\end{defi}

Let $\C$ be an additive category. It is said that $\C$ is {\bf Krull-Schmidt} if any object $X\in\T$ has a finite decomposition $X=\oplus_{i=1}^n\,X_i$ satisfying that each $X_i$ is indecomposable with local endomorphism ring $\End_\C(X_i).$ An idempotent $e=e^2\in\End_\C(X)$ 
{\bf splits} if there are morphism $u:X\to Y$ and $v:Y\to X$ satisfying $e=vu$ and $1_Y=uv.$
\

 The following result is well-known and a proof can be found, for example, in \cite{Chen}: An 
additive category $\C$ is Krull-Schmidt if and only if any idempotent in $\C$ splits and $\End_\C(X)$ is a semi-perfect ring for any $X\in \C.$ In this case, any object $X\in\C$ has a unique (up to order) finite direct decomposition $X=\oplus_{i=1}^n\,X_i$ satisfying that each $X_i$ is indecomposable with local endomorphism ring $\End_\C(X_i).$

\begin{defi} A category $\ltt{T}$ is said to be an {\bf artin triangulated
$R$-category} if the following conditions hold.
\begin{enumerate}
\item [(a)] $\ltt{T}$ is a triangulated $R$-category, where $R$ is
an artinian ring.
\item [(b)] $\ltt{T}$ is Hom-finite and Krull-Schmidt.
\end{enumerate}
\end{defi}

 Let $\Lambda$ be an artin $R$-algebra. It is also well-known 
that the bounded derived category $D^b(\Lambda),$ of complexes in $\modu\,(\Lambda),$ is an 
artin triangulated $R$-category (see, for example, in \cite[Theorem B.2]{Chen}).

\begin{pro}\label{ProAusl} Let $\ltt{T}$ be an artin triangulated
$R$-category, $A\in \ltt{T},$ $\Gamma:=\End_\T(A)^{op}$  and the evaluation functor at $A,$ 
$e_{A}:=\Hom_\T(A,-):\T \to \Modu\,(\Gamma).$ Then, the following conditions hold.
\begin{enumerate}
\item [(a)] $\Gamma$ is an artin $R$-algebra.
\item [(b)] The restriction, $e_{A}:\T \to \modu\,(\Gamma),$ is well 
defined and induces an equivalence of categories $\add\,(A)\stackrel{\sim}{\to}\proj\,(\Gamma).$
\item [(c)] $e_A:\Hom_\T(Z,X)\to \Hom_\Gamma(e_A(Z),e_A(X))$ is an isomorphism of $R$-modules for any $Z\in\add\,(A)$ and $X\in\T.$
\end{enumerate}
\end{pro}
\begin{dem} The proof done by M. Auslander (see \cite{Aus}) can be easily extended to the context of an artin triangulated $R$-category.
\end{dem}

\begin{lem}\label{triangulo-universal} Let $\ltt{T}$ be a Hom-finite triangulated $R$-category, 
 and let $A,C\in \T$ be such that
$\Hom_{\ltt{T}}(C[-1], A)\neq 0$. Then, the following conditions holds.
\begin{itemize}
\item [(a)] There exists a not splitting distinguished triangle in
$\ltt{T}$
$$\eta_{C,A}:\quad \xymatrix{A^{n}\ar[r]^{f} & E\ar[r]^{g} & C\ar[r]^{h} &
A^{n}[1]}$$ such that
$\Hom_{\ltt{T}}(-h[-1],A):\Hom_{\ltt{T}}(A^{n},A) \to \Hom_{\ltt{T}}(C[-1],A)$ is surjective, where
$n:=\ell_{R}(\Hom_{\ltt{T}}(C[-1],A))$.

\item [(b)] If $\Hom_{\ltt{T}}(A, A[1])=0$ then $\Hom_{\ltt{T}}(E, A[1])=0$.
\end{itemize}
\end{lem}
\begin{dem} (a) Since $n:=\ell_{R}(\Hom_{\ltt{T}}(C[-1],A)),$ it follows that there exists 
a family  $\{h_{i}\}_{i=1}^{n}$ of $R$-generators
in $\Hom_{\ltt{T}}(C,A[1])$. Hence, for each $i\in[1,n],$ we have the corresponding 
distinguished triangle
$$\eta_{i}:\xymatrix{A\ar[r]^{f_i} & B_{i}\ar[r]^{g_{i}} &
C\ar[r]^{h_{i}} & A[1]}.$$ By taking $\xi:=\oplus_{i=1}^{n}\,\eta_{i},$ we obtain the distinguished triangle 
$\xi: A^{n}\stackrel{}{\to}\oplus_{i=1}^{n}\,B_{i}\stackrel{}{\to} C^{n}\stackrel{}{\to} A^{n}[1].$ Let $\Delta:C
\fl C^{n}$ be the diagonal morphism. Then, by base change (see \ref{basecobase}), we get the
following commutative diagram  
$$\eta_{C,A}: \xymatrix{A^{n}\ar[r]^{f}\ar@{=}[d] & E\ar[r]^{g}\ar[d]^{K} & C\ar[r]^{h}\ar[d]^{\Delta} &
A^{n}[1]\ar@{=}[d]\\
A^{n}\ar[r] & \oplus_{i=1}^{n}\,B_{i}\ar[r] & C^{n}\ar[r] &
A^{n}[1],}$$ 
where the rows are distinguished triangles. Consider now, the following commutative diagram in
$\ltt{T}$
$$\xymatrix{C[-1]\ar[r]^{-h[-1]}\ar[d]_{\Delta[-1]} & A^{n}\ar[r]^{f}\ar@{=}[d] & E\ar[r]^{g}\ar[d]^{K}
 & C\ar[r]^{h}\ar[d]^{\Delta} &
A^{n}[1]\ar@{=}[d]\\
C^{n}[-1]\ar[r]^{-\varphi[-1]}\ar[d]_{\pi''_{i}[-1]}
& A^{n}\ar[r]\ar[d]^{\pi_{i}} &
\oplus_{i=1}^{n}\,B_{i}\ar[r]\ar[d]^{\pi'_{i}} & C^{n}\ar[r]^{\varphi}\ar[d]^{\pi''_{i}} & A^{n}[1]\ar[d]^{\pi_{i}[1]}\\
C[-1]\ar[r]^{-h_{i}[-1]} & A\ar[r] & B_{i}\ar[r]
& C\ar[r] & A[1],}$$ where $\pi_{i}$, $\pi_{i}'$ and
$\pi_{i}''$ are the corresponding canonical projections of the
direct sum. By the preceding diagram, we have that $\pi_{i}(-h[-1])=-h_{i}[-1]$ for all 
$i\in[1,n]$ and since the shift $[1]:\T\to \T$ is an 
$R$-functor, we get that the set $\{h_{i}[-1]\}_{i=1}^{n}$ is an $R$-generator of
$\Hom_{\ltt{T}}(C[-1],A)$. Thus the map
$$\Hom_{\ltt{T}}(-h[-1],A):\Hom_{\ltt{T}}(A^{n},A) \to \Hom_{\ltt{T}}(C[-1],A)$$ is surjective. Finally, by using the fact that $h_i\neq 0$ for 
each $i,$ it follows that $h\neq 0$ and therefore the triangle $\eta_{C,A}$ does not split.
\

(b) Let $\Hom_{\ltt{T}}(A, A[1])=0.$ Applying $\Hom_{\ltt{T}}(-,A)$ to the triangle
$\eta_{C,A},$ from the item (a), we have the following exact sequence
$$\xymatrix{(A^{n},A)\ar[rr]^{(-h[-1],A)} & &
(C[-1],A)\ar[r]   & (E[-1],A)\ar[r] &
(A^{n}[-1],A)}.$$ But, since
$\Hom_{\ltt{T}}(A^{n}[-1],A)=0$ and the map $\Hom_{\ltt{T}}(-h[-1],A)$ is surjective, it follows that
$\Hom_{\ltt{T}}(E[-1],A)=0$.
\end{dem}

\begin{pro}\label{triangulochido} Let $\ltt{T}$ be an artin triangulated
$R$-category and let $\eta:\,A\stackrel{s}{\to} B\stackrel{g}{\to} C\stackrel{\Psi}{\to} A[1]$
be a not splitting distinguished triangle such that
$\Hom_{\ltt{T}}(A,C)=\Hom_{\ltt{T}}(A[1],C)=0$ and $C$ is an
indecomposable object. Then, there exists a not splitting distinguished triangle
$A'\to B'\to  C\to  A'[1]$  such that $A'$ is a direct
summand of $A$ and $B'$ is an indecomposable direct summand of $B.$
\end{pro}
\begin{dem} Denote by $\alpha(B)$ the number of indecomposable direct summands
that appear in a decomposition of $B$ as direct sum of
indecomposables. The proof will be done by induction on
$\alpha(B)$. If $\alpha(B)=1$, there is nothing to prove.
\

Let $\alpha(B)>1$. Consider a decomposition $B=B_{1}\oplus B_{2}$
with $B_1$ indecomposable. Then, the triangle $\eta$ can be written as follows
$$\eta:\, \xymatrix{A\ar[rr]^{\left(\begin{array}{c}
s_1\\
s_2
\end{array}\right)} && B_{1}\oplus B_{2}\ar[rr]^{\left(\begin{array}{cc}
g_1 & g_2
\end{array}\right)} &&
 C\ar[rr] &&  A[1]}.$$
Applying $\Hom_{\ltt{T}}(-,C)$ to $\eta$, we have the following
exact sequence
$$\xymatrix{\Hom_{\ltt{T}}(A[1],C)\ar[r]  & \Hom_{\ltt{T}}(C,C)\ar[r]^{(g,C)} &
\Hom_{\ltt{T}}(B,C)\ar[r] & \Hom_{\ltt{T}}(A,C)}.$$ Since
$\Hom_{\ltt{T}}(A,C)=\Hom_{\ltt{T}}(A[1],C)=0$, it follows that
$\Hom_{\ltt{T}}(g,C)$ is an isomorphism. Let us consider the morphisms $(g_1,0):B\to C$
and $(0,g_2): B \to C$. Therefore there exists $f,\, f':C \to C$
such that $f(g_{1},g_{2})=(g_{1},0)$ and
$f'(g_{1},g_{2})=(0,g_{2})$. So, we get the following equalities
\begin{eqnarray*}
fg_{1}&=&g_{1}, \\
fg_{2}&=&0,     \\
f'g_{1}&=&0,    \\
f'g_{2}&=&g_{2}.
\end{eqnarray*}
Observe that $\Hom_{\ltt{T}}(g,C)(f+f')=(g_1,g_2),$ and since $\Hom_{\ltt{T}}(g,C)(1_C)=(g_1,g_2)$, we have that
$f+f'=1_C.$ We claim now that $f$ and $f'$ are idempotents. Indeed, we see, first, that $ff'=f'f=0$. The equality $ff'=0$ follows from the 
fact that  $\Hom_{\ltt{T}}(g,C)(ff')=(0,0)$ since $\Hom_{\ltt{T}}(g,C)$ is an isomorphism, and similarly we also get that $f'f=0$.\\
Now, from the equality $f+f'=1_C$, we get that
$f^{2}+ff'=f$ and then $f^{2}=f$. Analogously, it can be shown that
$f'^{2}=f'$. Furthermore, since $\ltt{T}$ is Krull-Schmidt and $C$ is
indecomposable, it follows that  either $f=0$ or $f'=0$. Hence, by the equalities listed above, we get that  either  $g_1=0$ or $g_2=0$.
\

Assume that $g_{1}=0$. Consider the following distinguished
triangles
$$\xymatrix{C[-1]\ar[r]^{h_2} & W'\ar[r]^{\delta'} &
B_{2}\ar[r]^{g_{2}} & C}\quad \text{and}  \quad \xymatrix{0\ar[r] &
B_1\ar[r]^{1_{B_1}} & B_1\ar[r] & 0},$$ 
where the first triangle is constructed by using the morphism $g_2.$ Thus, by taking their direct sum,
we get the following distinguished triangle
$$\xymatrix{C[-1]\ar[rr]^{\left(\begin{array}{c}
0\\
h_2
\end{array}\right)}  && B_{1}\oplus W'\ar[rr]^{\left(\begin{array}{cc}
1 & 0\\
0 & \delta'
\end{array}\right)}  & & B_{1}\oplus
B_2\ar[rr]^{\left(\begin{array}{cc} 0 & g_2
\end{array}\right)}  && C}.$$ So, we can construct the following
commutative diagram
$$\xymatrix{C[-1]\ar[rr]\ar@{=}[d] && A\ar[rr]^{\left(\begin{array}{c}
s_1\\
s_2
\end{array}\right)} && B_1\oplus B_2\ar[rr]^{\left(\begin{array}{cc}
0 & g_2
\end{array}\right)}\ar@{=}[d]  && C\ar@{=}[d]\\
C[-1]\ar[rr]^{\left(\begin{array}{c}
0\\
h_2
\end{array}\right)}  && B_{1}\oplus W'\ar[rr]^{\left(\begin{array}{cc}
1 & 0\\
0 & \delta'
\end{array}\right)}  & & B_{1}\oplus
B_2\ar[rr]^{\left(\begin{array}{cc} 0 & g_2
\end{array}\right)}  && C,}$$ 
where the rows are distinguished triangles. Hence, there exists an isomorphism $\xi:A \to
B_1\oplus W'$  inducing an isomorphism of triangles. In
particular, $W'$ is a direct summand of $A$. On the other hand, we
have the following commutative diagram
$$\xymatrix{W'\ar[r]^{\delta'}  & B_{2}\ar[r]^{g_2}\ar[d] & C\ar[r]^{-h_2[1]}\ar@{=}[d] &
 W'[1]\\
A\ar[r]^{s} & B\ar[r]^{g} & C\ar[r]^{\Psi} & A[1],}$$ 
where the rows are distinguished triangles. From
this diagram, we get a morphism $\beta':W' \fl A$ inducing a
morphism of triangles. Consider the following distinguished
triangle
$$\eta':\quad\xymatrix{W'\ar[r]^{\delta'} & B_2\ar[r]^{g_2} & C\ar[r]^{-h_2[1]} & W'[1].}$$ 
We assert that $\eta'$ does not split. Indeed, if $\eta'$ splits, we have that 
$-h_2[1]=0$  and then $\Psi=\beta'[1](-h_{2}[1])=0;$ thus the triangle $\eta$ splits, which is a contradiction proving that 
$\eta'$ does not split. Moreover $\Hom_{\ltt{T}}(W',C)=\Hom_{\ltt{T}}(W'[1],C)=0$ since
$W'$ is a direct summand of $A$. We also have that
$\alpha(\beta_{2})< \alpha(\beta)$. Hence, by induction, the result follows. Finally, the case $g_{2}=0$ is analogous. This completes the proof.
\end{dem}
\vspace{.4cm}

\section{filtered objects in a triangulated category}

Let $\ltt{X}$ be a class of objects in a triangulated category
$\ltt{T}$. It is said that an object $M\in\ltt{T}$ admits an
{\bf $\ltt{X}$-filtration} if there is a family of
distinguished triangles $\eta=\{\eta_{i}\;:\;M_{i-1}\to M_{i}\to X_{i}\to
 M_{i-1}[1]\}_{i=0}^{n}$ such that $M_{-1}=0=X_{0}$, $M_{n}=M$ and
$X_{i}\in \ltt{X}$ for $i\geq 1$. In such a case, it is defined the
lengths: $\ell_{\ltt{X},\eta}(M):=n$ and
$\ell_{\ltt{X}}(M):=\mini\{\ell_{\ltt{X},\eta}(M)\mid \eta
\,\,\text{is an}\,\,\ltt{X}\text{-filtration of}\, \,M\}$.
Finally, it is denoted by  $\F(\ltt{X})$ the class of objects
$M\in \ltt{T}$ for which there exists an $\ltt{X}$-filtration.

\begin{rk}\label{filtraex} For a triangulated category $\ltt{T}$ and
a class $\ltt{X}$ of objects in  $\ltt{T},$ the following statements hold.
\begin{enumerate}
\item [(a)] $\F(\ltt{X})=\cup_{n\in
\mathbb{N}}\;\F_{n}(\ltt{X})$, where
$\F_{0}(\ltt{X}):=\{0\}$ and
$\F_{n}(\ltt{X}):=\F_{n-1}(\ltt{X})\ast \ltt{X}$  for
$n\geq 1$.
\item [(b)] $\ell_{\ltt{X}}(M)=\mini\,\{n\in \N\mid M\in
\F_{n}(\ltt{X})\}$ for any $M\in \F(\ltt{X}).$

\item [(c)] $\F(\ltt{X}[i])=\F(\ltt{X})[i]$ for all $i\in \Z.$ Indeed, it can be
seen that  $(\ltt{X}\ast\ltt{Y})[i]=(\ltt{X}[i])\ast(\ltt{Y}[i])$ for any classes $\ltt{X}$ and  $\ltt{Y}$ of objects in $\ltt{T}$. Hence, (c) follows from (a).
\end{enumerate}
\end{rk}

\begin{lem}\label{extensionesclosed} Let $\ltt{X}$ be a class of objects in a triangulated 
category $\ltt{T}$. Then, the class $\F(\ltt{X})$ is closed under extensions.
\end{lem}
\begin{dem} Let $A\to B\to C\to A[1]$ be a
distinguished triangle in $\T$ with $A$ and $C$ in $\F(\ltt{X})$.
The proof will be done by induction on  $n:=\ell_{\ltt{X}}(C)$. If $C=0$, we have that 
$A\simeq B$ and hence $B\in\F(\ltt{X})$.\\
If $\ell_{\ltt{X}}(C)=1$ then $C\simeq X\in \ltt{X},$ and therefore
an $\ltt{X}$-filtration of $B$ can be done by adding the triangle $A\to B\to C\to A[1]$ to an 
$\ltt{X}$-filtration of $A$.
\

Suppose that $\ell_{\ltt{X}}(C)>1$. Consider a minimal $\ltt{X}$-filtration of $C$,
$$\{\eta_{i}: \xymatrix{C_{i-1}\ar[r] &
C_{i}\ar[r] & X_{i}\ar[r] & C_{i-1}[1]}\}_{i=0}^{n}.$$  By
base change (see \ref{basecobase}), we obtain the following commutative diagram in $\T$
$$\xymatrix{&  X_{n}[-1]\ar@{=}[r]\ar[d] & X_{n}[-1]\ar[d]\\
A\ar[r]\ar@{=}[d] & B_{n-1}\ar[r]\ar[d] & C_{n-1}\ar[r]\ar[d] & A[1]\ar@{=}[d]\\
A\ar[r] & B\ar[r]\ar[d] & C\ar[r]\ar[d] &  A[1]\\
& X_{n}\ar@{=}[r] & X_{n}}$$
where the rows and columns are distinguished triangles and $\ell_{\ltt{X}}(C_{n-1})<\ell_{\ltt{X}}(C)$. Applying induction to the first row of 
the preceding diagram, we get that $B_{n-1}\in \F(\ltt{X})$. Therefore,
an $\ltt{X}$-filtration of $B$ is given by adding the triangle
$B_{n-1}\to B\to X_{n}\to B_{n-1}[1]$ to an $\ltt{X}$-filtration of $B_{n-1}$.
\end{dem}

\begin{lem}\label{filtradogeneral} Let $\ltt{Y}$ and $\ltt{Z}$ be  classes of objects in a triangulated category
$\ltt{T}.$ If $\Hom_{\ltt{T}}(\ltt{Y},\ltt{Z}[i])=0$ for some $i\in\Z,$ then
$\Hom_{\ltt{T}}(\F(\ltt{Y}),\F(\ltt{Z})[i])=0$.
\end{lem}
\begin{dem} Let $\Hom_{\ltt{T}}(\ltt{Y},\ltt{Z}[i])=0$ for some $i\in\Z.$ By \ref{filtraex} (c), it is enough to prove the result only for the case 
$i=0.$ So, we assume that $\Hom_{\ltt{T}}(\ltt{Y},\ltt{Z})=0$ and we prove that $\Hom_{\ltt{T}}(\F(\ltt{Y}),\F(\ltt{Z}))=0.$
\

Let $N\in\F(\ltt{Y})$ and $M\in\F(\ltt{Z})$. We will show, by induction on
$\ell_{\ltt{Y}}(N),$ that $\Hom_{\ltt{T}}(N,M)=0$. In order to do that, we also can assume that $M\neq0$ and $N\neq 0$.\\
If $\ell_{\ltt{Y}}(N)=1$ then $N\simeq Y\in \ltt{Y}$ and so, by induction on $\ell_{\ltt{Z}}(M),$ it can be seen that 
$\Hom_{\ltt{T}}(N,M)=0$.\\
Suppose that $n:=\ell_{\ltt{Y}}(N)>1$. Then, there exists a distinguished triangle
$$\eta_{n}:\quad \xymatrix{N_{n-1}\ar[r] & N\ar[r] &
Y_{n}\ar[r] & N_{n-1}[1]}$$ such that $N_{n-1}\in
\F(\ltt{Y})$, $Y_{n}\in \ltt{Y}$ and
$\ell_{\ltt{Y}}(N_{n-1})=n-1$. Applying $\Hom_{\ltt{T}}(-,M)$ to
the triangle $\eta_{n}$, we get the exact sequence
$$\xymatrix{\Hom_{\ltt{T}}(Y_{n},M)\ar[r] &
\Hom_{\ltt{T}}(N,M)\ar[r] & \Hom_{\ltt{T}}(N_{n-1},M).}$$ By
induction, we have that $\Hom_{\ltt{T}}(N_{n-1},M)=0=\Hom_{\ltt{T}}(Y_{n},M)$, and therefore
$\Hom_{\ltt{T}}(N,M)=0$. 
\end{dem}

\begin{cor} \label{filtraper} Let $\ltt{X}$ be a class of objects in a
triangulated category $\ltt{T}.$ Then
$^{\perp}\ltt{X}=\,^{\perp}\F(\ltt{X})$.
\end{cor}
\begin{dem} It is enough to prove that ${}^{\perp} \ltt{X}\subseteq{}^{\perp}\F(\ltt{X}),$ since the other inclusion ${}^{\perp}\F(\ltt{X})
\subseteq {}^{\perp} \ltt{X}$ follows easily from the fact that $\ltt{X}\subseteq \F(\ltt{X}).$
\

Let $Y\in\, {}^{\perp}\ltt{X}$ and  $Z\in\F(\ltt{X})$. Then, by \ref{filtradogeneral}, it follows
that $\Hom_{\ltt{T}}(Y,Z)=0$, since $\Hom_{\ltt{T}}(Y,-)|_\ltt{X}=0.$ Thus $Y\in{}^{\perp}\F(\ltt{X})$ proving that 
$^{\perp}\ltt{X}=\,^{\perp}\F(\ltt{X})$.
\end{dem}

\begin{lem}\label{intercambio} Let $\ltt{T}$ be a triangulated category. If there are 
two distinguished triangles $Z\to Y\to \theta_{1}\to  Z[1]$ and $Y\to X\to \theta_{2}\to Y[1]$ such that 
$\Hom_{\ltt{T}}(\theta_{2},\theta_{1}[1])=0,$ then there exist two distinguished triangles as follows
$Z\to W\to \theta_{2}\to Z[1]$ and $W\to X\to \theta_{1}\to  W[1].$
\end{lem}
\begin{dem} Let $Z\to Y\to \theta_{1}\to  Z[1]$ and $Y\to X\to \theta_{2}\to Y[1]$ be 
distinguished triangles such that $\Hom_{\ltt{T}}(\theta_{2},\theta_{1}[1])=0.$
By co-base change (see \ref{basecobase}), we have the following commutative diagram
$$\xymatrix{& Z\ar@{=}[r]\ar[d] & Z\ar[d]\\
\theta_{2}[-1]\ar[r]\ar@{=}[d] & Y\ar[r]\ar[d] &
X\ar[r]\ar[d]
& \theta_{2}\ar@{=}[d]\\
\theta_{2}[-1]\ar[r] & \theta_{1}\ar[r]\ar[d] &
C\ar[r]\ar[d] & \theta_{2}\\
& Z[1]\ar@{=}[r] & Z[1],}$$ 
where the rows and columns are distinguished triangles. Using the fact that $\Hom_{\ltt{T}}(\theta_{2},\theta_{1}[1])=0,$ it follows that $\eta:\;\theta_{1}\to C\to \theta_{2}\to\theta_{1}[1]$ splits, and hence we get the following distinguished
triangle $\eta':\;\theta_{2}\to C\to \theta_{1}\to \theta_{2}[1]$. Then, by base change (see \ref{basecobase}), we
obtain the following commutative diagram 
$$\xymatrix{& \theta_{1}[-1]\ar@{=}[r]\ar[d] & \theta_{1}[-1]\ar[d]\\
Z\ar[r]\ar@{=}[d] & W\ar[r]\ar[d] & \theta_{2}\ar[r]\ar[d]
& Z[1]\ar@{=}[d]\\
Z\ar[r] & X\ar[r]\ar[d] &
C\ar[r]\ar[d] & Z[1]\\
& \theta_{1}\ar@{=}[r] & \theta_{1},}$$
where the rows and columns
are distinguished triangles. Hence, the required distinguished triangles are
$Z\to W\to \theta_{2}\to Z[1]$ and $W\to X\to \theta_{1}\to W[1].$
\end{dem}

\begin{lem}\label{filtraagrup} Let $\T$ be a triangulated category and $\theta\in \ltt{T}$ with $\Hom_{\ltt{T}}(\theta,\theta[1])=0,$ and let 
$\eta=\{\eta_{i}:\;M_{i-1}\to M_{i}\to \theta\to M_{i-1}[1]\}_{i=1}^{n}$ be a family of 
distinguished triangles. Then, for each $k\in [1,n],$ there exists a distinguished triangle 
$\xi_{k}:\; M_0\to M_{k}\to \theta^{k}\to M_{0}[1].$
\end{lem}
\begin{dem} We will proceed by induction on $k$. For $k=1,$ we have that 
$\xi_{1}:=\eta_{1}$ is the required triangle.
\

Let $k>1.$ Suppose we have $\xi_{k-1}$. By co-base change (see \ref{basecobase}), we get the 
following commutative diagram
$$\xymatrix{& M_{0}\ar@{=}[r]\ar[d] & M_{0}\ar[d]\\
\theta[-1]\ar[r]\ar@{=}[d] & M_{k-1}\ar[r]\ar[d] &
M_{k}\ar[r]\ar[d]
& \theta \ar@{=}[d]\\
\theta[-1]\ar[r] & \theta^{k-1}\ar[r]\ar[d] &
L_{k}\ar[r]\ar[d] & \theta\\
& M_{0}[0]\ar@{=}[r] & M_{0}[1],}$$ where the rows and columns are distinguished triangles. Since
$\Hom_{\ltt{T}}(\theta,\theta[1])=0$, the lower triangle
of the last diagram splits and hence $L_{k}\simeq
\theta^{k}$. Therefore, the second column of
the above diagram, is the required triangle $\xi_{k}.$
\end{dem}
\vspace{.4cm}

Let $\ltt{T}$ be a triangulated category and $\Theta=\{\Theta(i)\}_{i=1}^{t}$ be a family of
objects in $\T.$ For a given $\Theta$-filtration $\xi=\{\xi_k\;:\;M_{k-1}\to M_k\to X_k\to M_{k-1}[1]\}_{k=0}^n$ of $M\in\F(\Theta),$ we shall denote by 
$[M:\Theta(i)]_\xi$ the {\bf $\xi$-filtration multiplicity} of $\Theta(i)$ in $M.$ That is 
$[M:\Theta(i)]_\xi$ is the cardinal of the set $\{k\in[0,n]\;|\; X_k\simeq\Theta(i)\}.$ 
In general, the filtration multiplicity could be depending on a given $\Theta$-filtration. Observe that 
$\ell_{\Theta,\xi}(M)=\sum_{k=1}^t\;[M:\Theta(i)]_\xi.$

\begin{pro}\label{filtracionordenada} Let $\Theta=\{\Theta(i)\}_{i=1}^{t}$ be
a family of objects in a triangulated category  $\ltt{T}$, and let
$\leq$ be a linear order on $[1,t]$ such that
$\Hom_{\ltt{T}}(\Theta(j),\Theta(i)[1])=0$ for all $j\geq i$.
If $\xi$ is a $\Theta$-filtration of $M\in \F(\Theta)$, then there is a $\Theta$-filtration
$\eta$ of $M$ and a family $\Xi$ of distinguished triangles
satisfying the following conditions.
\begin{enumerate}
\item [(a)] $m(i):=[M:\Theta(i)]_{\xi}=[M:\Theta(i)]_{\eta}$ for all $i\in [1,t].$

\item [(b)] The family $\eta$ is ordered, that is,
$$\eta=\{\eta_{i}\;:\;M_{i-1}\to M_{i}\to \Theta(k_{i})\to M_{i-1}[1]\}_{i=0}^{n}$$ with
$\Theta(k_{0}):=0$, $M_{-1}:=0$ and $k_{n}\leq k_{n-1}\leq \cdots
\leq k_{1}$ in $([1,t],\leq)$.

\item [(c)] $\Xi=\{\Xi_{i}\,:\,M_{i-1}'\to M_{i}'\to \Theta(\lambda_{i})^{m(\lambda_{i})}\to 
M_{i-1}'[1]\}_{i=0}^{d},$  $\{\Theta(\lambda_{i})\}_{i=1}^{d}$
is the set consisting of the different $\Theta(j)$ appearing in
the $\Theta$-filtration $\xi$ of $M.$  Moreover
$\Theta(\lambda_{0}):=0,$ $M_{-1}':=0,$ $M'_d=M$ and $\lambda_{d}< \lambda_{d-1}<
\cdots < \lambda_{1}$ in $([1,t],\leq)$,
\end{enumerate}
\end{pro}
\begin{dem} Let $\xi$ be a $\Theta$-filtration of $M\in \F(\Theta).$ We can assume that 
$M\neq 0$ since the result is trivial in this case. 
\

We start by proving (a) and (b), proceeding by induction on
$n:=\ell_{\Theta,\xi}(M).$ If $n=1$, the $\Theta$-filtration $\xi$ is already ordered and hence
$\eta:=\xi$ satisfies the required properties. Let $n\geq 2$ and
$\xi:=\{\xi_{i}\;:\;M_{i-1}\to M_{i}\to \Theta(k_{i})\to M_{i-1}[1]\}_{i=0}^{n}$ be the 
$\Theta$-filtration of $M$. Since $\xi':=\xi-\{\xi_{n}\}$ is a
$\Theta$-filtration of $M_{n-1}$  and
$\ell_{\Theta,\xi'}(M_{n-1})=n-1$, by induction there is an ordered
$\Theta$-filtration $\eta'=\{\eta'_{i}\;:\;M'_{i-1}\to M'_{i}\to \Theta(k'_{i})\to M'_{i-1}[1]\}_{i=0}^{n-1}$ of $M_{n-1}$ with
$k'_{n-1}\leq k'_{n-2}\leq \cdots \leq k'_{1}$ and
$[M_{n-1}:\Theta(i)]_{\xi'}=[M_{n-1}:\Theta(i)]_{\eta'}$ $\;\forall\, i$. If $k_{n}\leq k'_{n-1}$ then $\eta:=\eta'\cup\{\xi_{n}\}$ satisfies the required conditions.\\
Suppose now that $k'_{n-1}<k_{n}.$ Let $l:=\maxi\{m\in [1,n-1]\;|\;
k'_{n-m}<k_{n}\}$. Observe that the $\Theta$-filtration
$\eta'\cup\{\xi_{n}\}$ is almost the one we want, the only triangle
that does not have its ordered multiplicity is precisely the
$\xi_{n}$. This can be rearranged by applying $l$-times
\ref{intercambio} to $\eta'\cup \{\xi_{n}\}$. 
\

(c) In order to construct $\Xi,$ we use the ordered $\Theta$-filtration $\eta$ from (b). We proceed as follows. For each $i\in[1,n]$, we group the $k_{i}$
that are the same and rename them by $\lambda_i$. So we get $\lambda_{d}< \lambda_{d-1}< \dots < \lambda_{1}$ on
$([1,t],\leq),$ and hence $\Theta(\lambda_{1}), \cdots, \Theta(\lambda_{d})$ are the different $\Theta(j)$ appearing in the $\Theta$-filtration $\eta$  of $M.$ Define
$s(i):=m(\lambda_{i})=[M:\Theta(\lambda_{i})],$ $\alpha(i):=\sum_{j=1}^{i}s(i)$ and $\alpha(0):=-1$.\\
We divide the filtration $\eta$ into the following pieces
$$\{\eta_{i}\;:\;\xymatrix{M_{i-1}\ar[r] &
M_{i}\ar[r] & \theta(\lambda_{l})\ar[r] & M_{i-1}[1]}\}_{i=\alpha(l-1)+1}^{\alpha(l)},$$ with $l\in[1,d]$. By
\ref{filtraagrup}, for each $l\in[1,d],$ we obtain the following distinguished triangle
$$\Xi_{l}:\quad \xymatrix{M_{\alpha(l-1)}\ar[r] & M_{\alpha(l)}\ar[r] &
\theta(\lambda_{l})^{s(l)}\ar[r] & M_{\alpha(l-1)[1]}.}$$ Hence, by setting  $\Xi_{0}:=\eta_{0}$ 
and $M_{i}':=M_{\alpha(i-1)}$ for $i\in[1, d],$ we conclude that the filtration
$\Xi=\{\Xi_{i}\}_{i=0}^{d}$ satisfies the required properties.
\end{dem}
\vspace{.4cm}

Let $\Theta=\{\Theta(i)\}_{i=1}^{t}$ be a family of objects in a triangulated category 
$\ltt{T}.$  We denote by $\Theta^\oplus$ the subcategory of $\T,$ whose objects are the 
finite direct sums of copies of objects in $\Theta.$

\begin{lem}\label{sumafiltra} Let $\Theta=\{\Theta(i)\}_{i=1}^{t}$ be a family of objects in a 
triangulated category $\ltt{T}.$ Then, the following statements hold.
\begin{itemize}
\item[(a)] $\F(\Theta)=\F(\Theta^\oplus).$
\item[(b)] If $\T$ is an artin triangulated $R$-category, then $\Theta^\oplus$ is functorially finite.
\end{itemize}
\end{lem}
\begin{dem} (a) Since $\Theta\subseteq \Theta^{\oplus}$, it follows that
$\F(\Theta)\subseteq\F(\Theta^{\oplus})$.
\

Let $M\in \F(\Theta^{\oplus}).$ We prove, by induction on
$m:=\ell_{\Theta^{\oplus}}(M),$ that $M\in\F(\Theta)$. If $m=1$,
then $M\in \Theta^{\oplus}$ and hence $M=\oplus_{i=1}^n\,\Theta(k_{i})^{m_i}.$ 
Since $\F(\Theta)$ is closed under extensions (see \ref{extensionesclosed}) and $\Theta(k_{i})\in
\F(\Theta)$,  we get that $M\in\F(\Theta)$.\\
Let $m>1.$ Then, there is a distinguished triangle
$M_{m-1}\to M\to \Theta(k_{m})^{\lambda(m)}\to M_{m-1}[1]$ with
$\ell_{\Theta^{\oplus}}(M_{m-1})=m-1.$ Hence, by induction, we get that $M_{m-1}\in
\F(\Theta)$. Therefore $M\in \F(\Theta)$ since
$\F(\Theta)$ is closed under extensions.
\

(b) The proof given in \cite[Proposition 4.2]{AusSma} can be easily extended to the context of an artin triangulated $R$-category.
\end{dem}

\begin{lem}\label{encajeFil} Let $\ltt{X}$ be a class of objects in a triangulated category
$\ltt{T}$ such that $0\in \ltt{X}$ and $\ltt{X}$ is closed under
isomorphisms. Then $\F_{n}(\ltt{X})=\ast_{i=1}^n\,\ltt{X}$ for $n\geq 1,$ and
$\F_{k}(\ltt{X})\subseteq \F_{k+1}(\ltt{X})$ for any $k\in\N.$
\end{lem}
\begin{dem} We have that $\F_{0}(\ltt{X}):=\{0\}$ and $\ltt{X}\subseteq
\F_{1}(\ltt{X})=\{0\}\ast\ltt{X}$. On the other hand, since
$\ltt{X}$ is closed under isomorphism, then
$\{0\}\ast\ltt{X}\subseteq \ltt{X}$. Hence,
$\F_{1}(\ltt{X})=\ltt{X}$ and so $\F_{2}(\ltt{X})=\ltt{X}\ast\ltt{X}$. Continuing in the same
way, we get that $\F_{n}(\ltt{X})=\ast_{i=1}^n\,\ltt{X}$ for $n\geq 1.$ Therefore, using the fact that the operation $\ast$ is associative, it follows that 
$\F_{k+1}(\ltt{X})=\ltt{X}\ast\F_{k}(\ltt{X}).$ Since $0\in \ltt{X},$ we conclude that
$\F_{k}(\ltt{X})\subseteq \F_{k+1}(\ltt{X})$ for any $k\in\N.$
\end{dem}
\vspace{.4cm}

The following result is a generalization, for triangulated categories, of the Ringel's result \cite[Theorem 1]{R}. The proof we give 
here uses the triangulated version of Gentle-Todorov's theorem due to Xiao-Wu Chen  
\cite[Theorem 1.3]{Chen2}.

\begin{teo}\label{TeoRingel} Let $\Theta=\{\Theta(i)\}_{i=1}^{t}$ be
a family of objects in an artin triangulated $R$-category  $\ltt{T}$, and let
$\leq$ be a linear order on the set $[1,t]$ such that\\
$\Hom_{\ltt{T}}(\Theta(j),\Theta(i)[1])=0$ for all $j\geq i$. Then
$\F(\Theta)=\ast_{i=1}^t\,\Theta^{\oplus}$ and it is functorially finite.
\end{teo}
\begin{dem} Let $\ltt{X}:=\Theta^{\oplus}.$ We assert that $\F_{t}(\ltt{X})=\F(\Theta)$. 
Indeed, by \ref{sumafiltra}, we have that
$\F(\ltt{X})=\F(\Theta)$ and hence
$\F(\Theta)=\cup_{n\in
\mathbb{N}}\,\F_{n}(\ltt{X})$.\\
Let $M\in \F(\Theta),$ and consider a $\Theta$-filtration $\xi$ of $M.$ Then, by \ref{filtracionordenada}(c), there is a family of distinguished triangles
$$\Xi=\{\Xi_{i}:\xymatrix{M_{i-1}'\ar[r] & M_{i}'\ar[r]
& \Theta(\lambda_{i})^{m(\lambda_{i})}\ar[r] & M_{i-1}'[1]}\}_{i=0}^{d},$$ where $\{\Theta(\lambda_{i})\}_{i=1}^{d}$
is the set of the different $\Theta(j)$ appearing in the $\Theta$-filtration 
$\xi$ of $M$, $\lambda_{d}<
\lambda_{d-1}< \cdots < \lambda_{1}$ and $M_{d}'=M$. Therefore $M\in
\F_{d}(\ltt{X})$ with $d \leq t$. Since
$\ltt{X}$ is closed under isomorphisms and contain the zero object,  by \ref{encajeFil}, 
it follows that $\F_{d}(\ltt{X})\subseteq
\F_{t}(\ltt{X})$. Thus 
$\F(\Theta)\subseteq \F_{t}(\ltt{X})$ and hence
$\F_{t}(\ltt{X})=\F(\Theta),$ proving our assertion.\\ 
By \ref{sumafiltra} (b), we know that $\ltt{X}$ is functorially finite. Furthermore, from \ref{encajeFil} and the assertion above, it follows
that $\F(\Theta)=\ast_{i=1}^t\,\ltt{X}.$ Hence the result follows from \cite[Theorem 1.3]{Chen2} and its dual.
\end{dem}
\vspace{.4cm}

\begin{defi} Let $\Theta=\{\Theta(i)\}_{i=1}^{n}$ be a family of objects in a triangulated category $\T.$ The {\bf $\Theta$-projective} objects in 
$\T$ is the class $\Pro(\Theta):={}^\perp\F(\Theta)[1].$ Dually, the {\bf $\Theta$-injective} objects in $\T$ is the 
class $\I(\Theta):=\F(\Theta)^\perp[-1].$ 
\end{defi}

Observe that, by \ref{filtraper} and its dual, we have that $\Pro(\Theta)={}^\perp\Theta[1]$  and $\I(\Theta)=\Theta^\perp[-1].$
\

In what follows, we use Ringel's ideas, in the paper \cite{R}, to proof that under certain conditions,
 $\Pro(\Theta)$ is a precovering class and $\I(\Theta)$ is a preenveloping one. To do that, we use the following two lemmas (compare 
with \cite[Lemma 3 and Lemma 4]{R}).

\begin{lem}\label{trianparcial} Let $\Theta=\{\Theta(i)\}_{i=1}^{n}$ be a
family of objects, in a  Hom-finite triangulated $R$-category $\T,$ such
that  $\Hom_{\ltt{T}}(\Theta(j),\Theta(i)[1])=0$ for $j\geq i$.
Consider $t\in [1,n]$ and $N\in \ltt{T}$  such that
$\Hom_{\ltt{T}}(\Theta(j), N[1])=0$ for $j>t$. Then, there
exists a distinguished triangle in $\T$
$$\xymatrix{N\ar[r] & N_{t}\ar[r] & \Theta(t)^{m}\ar[r] & N[1]},$$
where $m:=\ell_{R}\Hom_{\ltt{T}}(\Theta(t),N[1])$ and
$\Hom_{\ltt{T}}(\Theta(j), N_{t}[1])=0$ for $j\geq t$.
\end{lem}
\begin{dem} If $\Hom_{\ltt{T}}(\Theta(t),N[1])=0$, the distinguished triangle we are looking for 
is $N\stackrel{1}{\to} N\to 0\to N[1].$
\

Let $\Hom_{\ltt{T}}(\Theta(t),N[1])\neq 0.$ Then, by the dual of \ref{triangulo-universal}, there is a distinguished triangle
$$\eta:\xymatrix{N\ar[r] & N_{t}\ar[r] & \Theta(t)^{m}\ar[r]^{h} &
 N[1]}$$ such that the map 
$\Hom_{\ltt{T}}(\Theta(t),h):\Hom_{\ltt{T}}(\Theta(t),\Theta(t)^{m})
\to\Hom_{\ltt{T}}(\Theta(t),N[1])$ is surjective. Applying
$\Hom_{\ltt{T}}(\Theta(j),-)$  to $\eta,$ we get the following
exact sequence
$$\xymatrix{(\Theta(j),\Theta(t)^{m})\ar[r] & (\Theta(j),N[1])\ar[r] & 
(\Theta(j),N_{t}[1])\ar[r] &
(\Theta(j),\Theta(t)^{m}[1]).}$$ Since
$\Hom_\T(\Theta(j),\Theta(t)[1])=0$ for $j\geq t$ and
$\Hom_\T(\Theta(j), N[1])=0$ for $j>t$, it follows that
$\Hom_\T(\Theta(j), N_{t}[1])=0$ for $j>t$. For $j=t$, we know that
$\Hom_{\ltt{T}}(\Theta(t),h)$ is an epimorphism and hence
$\Hom_\T(\Theta(t), N_{t}[1])=0;$ proving the lemma.
\end{dem}

\begin{lem}\label{triancasi} Let $\Theta=\{\Theta(i)\}_{i=1}^{n}$ be a
family of objects, in a Hom-finite triangulated $R$-category $\T$, such
that $\Hom_{\ltt{T}}(\Theta(j),\Theta(i)[1])=0$ for $j\geq i$.
Consider $t\in [1,n]$ and $N\in \ltt{T}$ such that
$\Hom_{\ltt{T}}(\Theta(j), N[1])=0$ for $j>t$. Then, there exists
a distinguished triangle in $\T$
$$\xymatrix{N\ar[r] & Y\ar[r] & X\ar[r] & N[1]}$$ with $X\in
\F(\{\Theta(i)\;|\; i\in [1,t]\})$ and $Y\in\I(\Theta).$
\end{lem}
\begin{dem} Since $\Hom_{\ltt{T}}(\Theta(j),N[1])=0$ for $j>t$, it follows from 
\ref{trianparcial} the existence of a distinguished triangle
$$\eta_{t+1}:\xymatrix{N\ar[r]^{\mu_{t}} & N_{t}\ar[r] & Q_{t}\ar[r] & N[1]}$$
with $Q_{t}:=\Theta(t)^{m_{t}}$ and $\Hom_{\ltt{T}}(\Theta(j),N_{t}[1])=0$ for $j\geq t$. 
Similarly, there is a distinguished triangle
$$\eta_{t}:\xymatrix{N_{t}\ar[r]^{\mu_{t-1}} & N_{t-1}\ar[r] & Q_{t-1}\ar[r] & N_{t}[1]}$$
with $Q_{t-1}:=\Theta(t-1)^{m_{t-1}}$ and $\Hom_{\ltt{T}}(\Theta(j),N_{t-1}[1])=0$ for $j\geq t-1$. Continuing this procedure, we
get distinguished triangles
$$\eta_{i}:\xymatrix{N_{i}\ar[r]^{\mu_{i-1}} & N_{i-1}\ar[r] & Q_{i-1}\ar[r] & N_{i}[1]}$$
with $Q_{i-1}:=\Theta(i-1)^{m_{i-1}}$ and $\Hom_{\ltt{T}}(\Theta(j), N_{i}[1])=0$ for $j\geq i$. In what follows, for
$\alpha_{r}:=\mu_{t-r}\dots\mu_{t-1}\mu_{t}$ with $0\leq r\leq
t-1$, we will construct, inductively, distinguished triangles
$$\xi_{r}: \xymatrix{N\ar[r]^{\alpha_{r}} & N_{t-r}\ar[r] &
X_{t-r}\ar[r] & N[1]}$$ with $X_{t-r}\in
\F(\{\Theta(i)\;|\; i\in [t-r,t]\})$ and $\Hom_{\ltt{T}}(\Theta(j),N_{t-r}[1])=0$ for 
$j\geq t-r$. If $r=0$, we set $\xi_{0}:=\eta_{t+1}$. Suppose that $r>0$ and that the triangle $\xi_{r}$ is already constructed. Consider the following diagram of co-base change (see \ref{basecobase})
$$\xymatrix{& \Theta(t-r-1)^{m_{t-r-1}}[-1]\ar@{=}[r]\ar[d] &  \Theta(t-r-1)^{m_{t-r-1}}[-1]\ar[d]\\
N\ar[r]^{\alpha_{r}}\ar@{=}[d] & N_{t-r}\ar[r]\ar[d]^{\mu_{t-r-1}}
& X_{t-r}\ar[r]\ar[d]
& N[1]\ar@{=}[d]\\
N\ar[r]^{\alpha_{r+1}} & N_{t-r-1}\ar[r]\ar[d] &
X_{t-r-1}\ar[r]\ar[d] & N[1] \\
& \Theta(t-r-1)^{m_{t-r-1}}\ar@{=}[r] &
\Theta(t-r-1)^{m_{t-r-1}}.}$$ By induction, we have that
$X_{t-r}\in \F(\{\Theta(i)\;|\;i\in [t-r,t]\})$. Thus,
$\Theta(t-r-1)^{m_{t-r-1}}$, $X_{t-r}\in \F(\{\Theta(i)\;|\;
i\in [t-r-1,t]\})$. Since $\F(\{\Theta(i)\;|\; i\in
[t-r-1,t]\})$ is closed under extensions, it follows that
$X_{t-r-1}\in \F(\{\Theta(i)\;|\; i\in [t-r-1,t]\})$. Moreover
$\Hom_{\ltt{T}}(\Theta(j),N_{t-r-1}[1])=0$ for $j\geq t-r-1$.
Therefore $\xi_{r+1}$ is the triangle from the second row of the
last diagram. Then, the required triangle is $\xi_{t-1}$.
\end{dem}

\begin{teo}\label{TeoRingel2} Let $\Theta=\{\Theta(i)\}_{i=1}^{n}$ be
a family of objects in an artin triangulated $R$-category  $\ltt{T}$, and let
$\leq$ be a linear order on the set $[1,t]$ such that\\
$\Hom_{\ltt{T}}(\Theta(j),\Theta(i)[1])=0$ for all $j\geq i$. Then, the following statements holds.
\begin{itemize}
 \item[(a)] For any object $X\in\T$ there are two distinguished triangles in $\T$
$$\begin{matrix} \xymatrix{X\ar[r] & Y_X\ar[r] & C_X\ar[r] &  X[1]} & \text{ with } Y_X\in\I(\Theta),\,C_X\in\F(\Theta),\\ {}\\
  \xymatrix{X[-1]\ar[r] & K_X\ar[r] & Q_X\ar[r] & X }     & \text{ with } Q_X\in\Pro(\Theta),\,K_X\in\F(\Theta).  
  \end{matrix}$$ 
\item[(b)] $\Pro(\Theta)$ is a precovering class and $\I(\Theta)$ is a preenveloping one in $\T.$
\end{itemize}
\end{teo}
\begin{dem} (a) For simplicity, we assume that the linear order $\leq$ on the set $[1,t]$ is the natural one. Furthermore, we only prove the existence 
of the first triangle, since the existence of the other one follows by duality. Let $X\in\T$ and $t:=n.$ Then, from \ref{triancasi}, we get a distinguished 
triangle $X\to Y_X\to C_X\to X[1]$ in $\T$ such that $Y_X\in\I(\Theta)$ and $C_X\in\F(\Theta).$
\

(b) We start proving that  $\I(\Theta)$ is a preenveloping class in $\T.$ Indeed, let $X\in \ltt{T}.$  Then, by (a), there is a
distinguished triangle
$$\xymatrix{X\ar[r]^{\beta} & Y_X\ar[r] & C_X\ar[r] & X[1]}$$ with $Y_X\in\I(\Theta)$ and $C_X\in\F(\Theta).$ We claim
that $\beta$ is an $\I(\Theta)$-preenvelope of $X$. To see that, we consider a morphism $\beta':X \to  Y'$ with $Y'\in\I(\Theta).$ Then, 
by co-base change (see \ref{basecobase}), we have the following
commutative diagram in $\T$
$$\xymatrix{C_X[-1]\ar[r]\ar@{=}[d] & X\ar[d]^{\beta'}\ar[r]^{\beta} & Y_X\ar[r]\ar[d]^{\gamma} &
C_X\ar@{=}[d]\\
C_X[-1]\ar[r]^{\alpha} & Y'\ar[r]^{u} & L\ar[r] & C_X,}$$
where the rows are distinguished triangles. Since  $Y'\in\I(\Theta)=\F(\Theta)^{\perp}[-1]$ and $C_X\in \F(\Theta),$ we get that
$\Hom_{\ltt{T}}(C_X[-1],Y')=0.$  Therefore $\alpha=0$ and thus $\beta'$ factors through $\beta;$ proving that 
$\beta$ is an $\I(\Theta)$-preenvelope of $X.$ Finally, the proof that $\Pro(\Theta)$ is a precovering class in $\T$ is rather similarly by 
using the second triangle in (a). 
\end{dem}

\section{Homological systems}

In this section, we introduce several homological systems of objects in a triangulated category $\T,$ over a linearly ordered finite set. This 
homological systems generalize the notion of stratifying systems (see \cite{ES,MMS1, MMS2,MSX}) in a module category. We recall that 
$[1,n]:=\{1,2,\cdots,n\}$ for any $n\in\Z^+.$

\begin{defi}\label{ss} A {\bf $\Theta$-system} $(\Theta,\leq)$ of size t, in a triangulated
category  $\ltt{T}\!,$ consists of the following data.
\begin{enumerate}
\item[(S1)] $\leq$ is a linear order on $[1,t].$
\item[(S2)] $\Theta=\{\Theta(i)\}_{i=1}^{t}$ is a family of indecomposable objects in $\T$.
\item[(S3)] $\Hom_{\ltt{T}}(\Theta(j),\Theta(i))=0$ for $j>i$.
\item[(S4)] $\Hom_{\ltt{T}}(\Theta(j),\Theta(i)[1])=0$ for $j\geq i$.
\item[(S5)] $\Hom_{\ltt{T}}(\Theta,\Theta[-1])=0$.
\end{enumerate}
\end{defi}

\begin{defi} \label{epss} A {\bf $\Theta$-projective system} $(\Theta,\Qb,\leq)$ of size t, in a triangulated
category  $\ltt{T}\!,$ consists of the following data.
\begin{enumerate}
\item[(PS1)] $\leq$ is a linear order on $[1,t].$
\item[(PS2)] $\Theta=\{\Theta(i)\}_{i=1}^{t}$ is a family of non-zero objects in $\ltt{T}.$ 
\item[(PS3)] $\Hom_{\ltt{T}}(\Theta(j),\Theta(i))=0$ for $j>i$.
\item[(PS4)] $\Qb=\{Q(i)\}_{i=1}^{t}$ is a family of indecomposable objects in  $\ltt{T}$ such 
that $Q:=\bigoplus_{i=1}^{t}\,Q(i)\in\,
{}^{\perp}\Theta[-1]\cap \,{}^{\perp}\Theta[1].$
\item [(PS5)] For every $i\in [1,t]$, there exists a distinguished
triangle in $\T$
$$\eta_{i}\;:\; \xymatrix{K(i)\ar[r] & Q(i)\ar[r]^{\beta_{i}} &
\Theta(i)\ar[r] & K(i)[1]}$$ such that $K(i)\in
\F(\{\Theta(j)\;|\; j>i\})$ and $\Hom_{\ltt{T}}(K(i)[1],\Theta(i))=0$.
\end{enumerate}
\end{defi}

\begin{defi} \label{eiss} A {\bf $\Theta$-injective system} $(\Theta,\Yb,\leq)$ of size t, in a triangulated
category  $\ltt{T}\!,$ consists of the following data.
\begin{enumerate}
\item[(IS1)] $\leq$ is a linear order on $[1,t].$
\item[(IS2)] $\Theta=\{\Theta(i)\}_{i=1}^{t}$ is a family of non-zero objects in $\ltt{T}.$ 
\item[(IS3)] $\Hom_{\ltt{T}}(\Theta(j),\Theta(i))=0$ for $j>i$.
\item[(IS4)] $\Yb=\{Y(i)\}_{i=1}^{t}$ is a family of indecomposable objects in  $\ltt{T}$ such 
that $Y:=\bigoplus_{i=1}^{t}\,Y(i)\in\,
\Theta^{\perp}[-1]\cap \Theta^{\perp}[1].$
\item [(IS5)] For every $i\in [1,t]$, there exists a distinguished
triangle in $\T$
$$\xi_{i}\;:\; \xymatrix{Z(i)[-1]\ar[r] & \Theta(i)\ar[r]^{\alpha_{i}} &
Y(i)\ar[r] & Z(i)}$$ such that $Z(i)\in
\F(\{\Theta(j)\;|\; j<i\})$ and $\Hom_{\ltt{T}}(\Theta(i),Z(i)[-1])=0$.
\end{enumerate}
\end{defi}

\begin{rk}\label{InjVsProj} A triple $(\Theta,\Yb,\leq)$ is a $\Theta$-injective system of size $t$, in a triangulated category $\T,$ 
if and only if $(\Theta^{op},\Yb^{op},\leq^{op})$ is a $\Theta^{op}$-projective system of size $t$ in the opposite triangulated 
category $\T^{op},$ where $\leq^{op}$ is the opposite order of $\leq$ in $[1,t].$ Therefore, any obtained result for  $\Theta$-projective systems
can  be transfered to the  $\Theta$-injective systems, and so, we could be dealing only with  $\Theta$-projective systems.
\end{rk}

\begin{pro}\label{anulasuperior} Let $(\Theta,\Qb,\leq )$ be a $\Theta$-projective system
of size $t,$ in a triangulated category $\ltt{T}$. Then, the
following conditions hold.
\begin{enumerate}
\item [(a)]$\Hom_{\ltt{T}}(K(j),\Theta(i))=0=\Hom_{\ltt{T}}(\Theta(j),\Theta(i)[1])$
for all $j\geq i$.
\item [(b)] $\Hom_{\ltt{T}}(\beta_{j},\Theta(i)):\Hom_{\ltt{T}}(\Theta(j),\Theta(i))\to
\Hom_{\ltt{T}}(Q(j),\Theta(i))$ is an isomorphism of abelian groups, for all $j\geq i$.
\item [(c)] If $\Hom_{\ltt{T}}(K(j)[2],\Theta(i))=0$ $\;\forall\, i,j\in [1,t]$, then
$\Hom_{\ltt{T}}(\Theta,\Theta[-1])=0$.
\end{enumerate}
\end{pro}
\begin{dem} (a) Let $j\geq i.$ Using the fact that $K(j)\in\F(\{\Theta(\lambda)\;|\;\lambda >j\})$ 
 and since
$\Hom_{\ltt{T}}(\Theta(\lambda),\Theta(i))=0$ for $\lambda >j\geq
i,$ it follows from \ref{filtradogeneral} that $\Hom_{\ltt{T}}(K(j),\Theta(i))=0.$
Consider the distinguished triangle given in \ref{epss} (PS5)
$$\eta_{j}\;:\;
\xymatrix{K(j)\ar[r] & Q(j)\ar[r]^{\beta_j} & \Theta(j)\ar[r] & K(j)[1]}.$$ Applying 
$\Hom_{\ltt{T}}(-,\Theta(i)[1])$ to $\eta_{j}$, we get the exact sequence
$$\xymatrix{(K(j)[1],\Theta(i)[1])\ar[r] &
(\Theta(j),\Theta(i)[1])\ar[r] & (Q(j),\Theta(i)[1])}.$$
 Thus, since $\Qb\subseteq{}^{\perp}\Theta[1]$ and $\Hom_{\ltt{T}}(K(j),\Theta(i))=0,$ it 
 follows from the sequence above that $\Hom_{\ltt{T}}(\Theta(j),\Theta(i)[1])=0$ for $j\geq i$.
\

(b) Let $j\geq i.$ Applying $\Hom_{\ltt{T}}(-,\Theta(i))$ to the above distinguished triangle
$\eta_{j}$, we get the exact sequence
$$\xymatrix{(K(j)[1],\Theta(i))\ar[r] &
(\Theta(j),\Theta(i))\ar[r]^{(\beta_{j},\Theta(i))} &
(Q(j),\Theta(i))\ar[r] & (K(j),\Theta(i)).}$$ 
We have that $\Hom_\T(\beta_{j},\Theta(i))$ is an epimorphism, since by (a) we know that
$\Hom_\T(K(j),\Theta(i))=0$ for $j\geq i.$ Since $\Hom_\T(K(i)[1],\Theta(i))=0$ (see \ref{epss} 
(PS5)), we conclude that $\Hom_\T(\beta_{i},\Theta(i))$ is an isomorphism.\\
Assume that $j>i$. Then $\Ker\,(\Hom_\T(\beta_{j},\Theta(i)))\subseteq \Hom_\T(\Theta(j),\Theta(i))=0$ and hence $\Hom_\T(\beta_{j},\Theta(i))$ is also an isomorphism.
\

(c) Let $i,j\in[1,t].$ Applying $\Hom_{\ltt{T}}(-,\Theta(i)[-1])$ to the above distinguished triangle
$\eta_{j}$, we get the exact sequence
$$\xymatrix{(K(j)[1],\Theta(i)[-1])\ar[r] &
(\Theta(j),\Theta(i)[-1])\ar[r] &
(Q(j),\Theta(i)[-1])}.$$ 
Using the fact that $\Qb\subseteq {}^\perp\Theta[-1]$ and since $\Hom_{\ltt{T}}
(K(j)[2],\Theta(i))=0,$ it follows that
$\Hom_{\ltt{T}}(\Theta(j),\Theta(i)[-1])=0;$ proving that $\Hom_{\ltt{T}}(\Theta,\Theta[-1])=0.$
\end{dem}

\begin{pro}\label{Basico1Proj}
Let $(\Theta,\Qb,\leq )$ be a $\Theta$-projective system of size $t,$ in an
artin triangulated $R$-category $\ltt{T}$. Then, the following statements hold.
\begin{enumerate}
\item [(a)] For each $i\in [1,t]$, the morphism $\beta_{i}:Q(i)\to \Theta(i),$ appearing in 
the triangle $\eta_i$ from \ref{epss} (PS5), is a $\Pro(\Theta)$-cover of $\Theta(i)$.
\item [(b)] Let $(\Theta,\Qb',\leq )$ be another $\Theta$-projective system of size 
$t,$ in $\ltt{T}.$ Then $\Qb'\simeq\Qb;$ that is, for each $i\in [1,t],$ there is an
isomorphism $\rho_{i}:Q(i) \to Q'(i)$ such that the following diagram in $\T$ commutes 
$$\xymatrix{Q(i)\ar[rd]_{\beta_{i}}\ar[rr]^{\rho_{i}} & & Q'(i)\ar[dl]^{\beta_{i}'}\\
 &  \Theta(i). & }$$
\end{enumerate}
\end{pro}
\begin{dem} (a) Let $i\in [1,t].$ We start by proving that $\beta_{i}:Q(i)\to \Theta(i)$ is right 
minimal. Firstly, we assert that $\beta_i\neq 0.$ Indeed, by \ref{anulasuperior} (b), we have 
that $\Hom_{\ltt{T}}(\beta_{i},\Theta(i)):\End_\T(\Theta(i))\to \Hom_\T(Q(i),\Theta(i))$ is 
an isomorphism. Thus $\beta_{i}=\Hom_{\ltt{T}}(\beta_{i},\Theta(i))(1_{\Theta(i)})\neq
0$ since $1_{\Theta(i)}\neq 0.$ Let $f: Q(i) \to Q(i)$ be such that  $\beta_i f=\beta_i$. Then 
$\beta_i=\beta_i f^{n}$ $\;\forall \, n\in \mathbb{N}^{+}$. Since $\beta_i\neq 0$, it follows 
that $f^{n}\neq 0$ $\forall\, n\in \mathbb{N}^{+}$. Using the fact that $Q(i)$ is 
indecomposable, we get  from \ref{ProAusl} (a), that $\End_\T(Q(i))$ is a local artin
R-algebra. Thus, $\rad\,(\End_\T(Q(i)))$ is nilpotent and coincides with the set of 
non-invertible elements of $\End_\T(Q(i))$. Since $f^{n}\neq 0$ $\;\forall\,n\in \mathbb{N}^{+}$, 
we conclude that $f\notin \rad\,(\End_\T(Q(i)))$ and therefore $f$ is invertible; proving that 
$\beta_{i}:Q(i)\to \Theta(i)$ is right minimal.\\
Finally, we prove that $\beta_{i}:Q(i)\to \Theta(i)$ is a $\Pro(\Theta)$-precover of $\Theta(i)$.
Let $g: X \to\Theta(i)$ be in $\T,$ with $X \in\Pro(\Theta)$. Applying $\Hom_{\ltt{T}}(X,-)$  to
the distinguished triangle $\eta_{i}$ from \ref{epss} (PS5), we get the exact sequence $$\xymatrix{(X,K(i))\ar[r]
&(X,Q(i))\ar[r]^{(X,\beta_{i})} & (X,\Theta(i))\ar[r] & (X,K(i)[1]).}$$ Since 
$X \in\Pro(\Theta)$ and $K(i)\in\F(\Theta)$, we conclude that $\Hom_\T(X,K(i)[1])=0;$ proving 
that $g$ factorizes through $\beta_{i},$ and thus $\beta_{i}$  is a $\Pro(\Theta)$-precover of $\Theta(i).$
\

(b) It is immediate from (a)
\end{dem}
\vspace{.4cm}

Let $(\Theta,\leq)$ be $\Theta$-system in a triangulated category  $\ltt{T}.$ A natural question here, is to ask for the existence of a family
$\Qb$ of objects in $\T$ such that $(\Theta,\Qb,\leq)$ is a $\Theta$-projective system. In order to do that, we will need the 
following results. Recall that, for any $a,b\in\Z$ with $a\leq b,$ we set $[a,b]:=\{x\in\Z\;|\;a\leq x\leq b\}.$

\begin{lem}\label{filtradocero} Let $(\Theta,\leq)$ be a $\Theta$-system of size $t,$ in a triangulated category $\ltt{T}$, where $\leq $ 
is the natural order on the set $[1,t]$. Then, the following statements hold.
\begin{enumerate}
\item [(a)] If $M\in \F(\{\Theta(j)\;|\; j\in[i,i+k]\})$,
$N\in \F(\{\Theta(r)\;|\; r>i+k\})$ and $L\in \F(\{\Theta(s)\;|\; s<i\})$, then $\Hom_{\ltt{T}}(N,M)=0$ and $\Hom_{\ltt{T}}(M,L)=0$.
\item [(b)] If $M\in \F(\{\Theta(j)\;|\;j\in[i,i+k]\})$,
$N\in \F(\{\Theta(r)\;|\; r\geq i+k\})$ and $L\in\F(\{\Theta(s)\;|\; s\leq i\})$, then $\Hom_{\ltt{T}}(N,M[1])=0$ and $\Hom_{\ltt{T}}(M,L[1])=0$.
\item [(c)] If $M, N\in \F(\Theta)$ then $\Hom_{\ltt{T}}(M,N[-1])=0$.
\end{enumerate}
\end{lem}
\begin{dem}
It follows immediately from \ref{filtradogeneral} and the
definition of stratifying system.
\end{dem}

\begin{pro}\label{existriangulos} Let $(\Theta,\leq )$ be a $\Theta$-system of size $t,$ in 
an artin triangulated $R$-category $\ltt{T}\!,$ and let $\leq$ be the natural order on $[1,t]$, $t> 1$ and $i\in[1,t]$. Then, for each 
$k\in[1,t-i],$ there exists a distinguished triangle in $\T$
$$\xi_{k}\;:\; \xymatrix{V_{k}\ar[r] & U_{k}\ar[r] &
\Theta(i)\ar[r] & V_{k}[1]}$$ satisfying the following
conditions:
\begin{enumerate}
\item [(a)] $U_{k}$ is indecomposable,
\item [(b)] $V_{k} \in \F(\{\Theta(j)\;|\;i<j \leq i+k\}),$
\item [(c)] $\Hom_{\ltt{T}}(U_{k},\Theta(j)[1])=0$ for $j\in[i,i+k]$.
\end{enumerate}
\end{pro}
\begin{dem} We will proceed by induction on $k$.\\
Let $k=1$. By definition, we have that $\Hom_{\ltt{T}}(\Theta(i+1), \Theta(i))=0$. If
$\Hom_{\ltt{T}}(\Theta(i),\Theta(i+1)[1])=0,$ the desired triangle is the following 
$$\xymatrix{0\ar[r] & \Theta(i)\ar[r]^{1} & \Theta(i)\ar[r]  &
0.}$$ Suppose that $\Hom_{\ltt{T}}(\Theta(i),\Theta(i+1)[1])\neq
0$. Then, by \ref{triangulo-universal}, there exists a not splitting distinguished 
triangle in $\T$
$$\xi\;:\; \xymatrix{\Theta(i+1)^{n}\ar[r]  & E\ar[r] & \Theta(i)
\ar[r] & \Theta(i+1)^{n}};$$ and moreover, we have that
$\Hom_{\ltt{T}}(E,\Theta(i+1)[1])=0$. Applying the functor
$\Hom_{\ltt{T}}(-,\Theta(i)[1])$ to $\xi$, we get the exact
sequence
$$\xymatrix{\Hom_{\ltt{T}}(\Theta(i),\Theta(i)[1])\ar[r]
& \Hom_{\ltt{T}}(E,\Theta(i)[1])\ar[r] &
\Hom_{\ltt{T}}(\Theta(i+1)^{n},\Theta(i)[1]).}$$ Since
$\Hom_{\ltt{T}}(\Theta(i),\Theta(i)[1])=0=
\Hom_{\ltt{T}}(\Theta(i+1)^{n},\Theta(i)[1])=0$, we conclude
that $\Hom_{\ltt{T}}(E,\Theta(i)[1])=0$. Moreover, since $\Hom_{\ltt{T}}
(\Theta(i+1)^{n},\Theta(i))=0=\Hom_{\ltt{T}}(\Theta(i+1)^{n},\Theta(i)[-1])=0$, it follows by \ref{triangulochido} the existence of a distinguished triangle
$$\xi':\xymatrix{\Theta(i+1)^{m}\ar[r]  & U_{1}\ar[r] & \Theta(i)\ar[r] &
\Theta(i+1)^{m}[1]}$$ with $m\leq n$ and $U_{1}$ an
indecomposable direct summand of $E$. Thus, the distinguished triangle 
$\xi_{1}:=\xi'$ satisfies the required conditions. Suppose now
that there exists a distinguished triangle
$$\xi_{k}:\quad \xymatrix{V_{k}\ar[r] & U_{k}\ar[r] &
\Theta(i)\ar[r] & V_{k}[1]}$$ satisfying the above required properties.
We construct the distinguished triangle $\xi_{k+1},$ from $\xi_{k},$ as follows. If $\Hom_{\ltt{T}}(U_{k},\Theta(i+k+1)[1])=0$, the triangle
$\xi_{k+1}:=\xi_{k}$ is the desired one.
\

Suppose that $\Hom_{\ltt{T}}(U_{k},\Theta(i+k+1)[1])\neq 0$. Then, by 
\ref{triangulo-universal}, there exists a not splitting distinguished
triangle in $\T$
$$\eta:\quad \xymatrix{\Theta(i+k+1)^{a}\ar[r] & U\ar[r] & U_{k}\ar[r] &
\Theta(i+k+1)^{a}[1]},$$ and furthermore we have that 
$\Hom_{\ltt{T}}(U,\Theta(i+k+1)[1])=0.$ Applying the functor
$\Hom_{\ltt{T}}(-,\Theta(i+k+s)[1])$ to $\eta$, with $s\in[-k,0],$ we get the exact sequence
$$(U_{k},\Theta(i+k+s)[1])\to (U,\Theta(i+k+s)[1])\to (\Theta(i+k+1)^{a},\Theta(i+k+s)[1]).$$ Since 
$\Hom_{\ltt{T}}(\Theta(i+k+1)^{a},\Theta(i+k+s)[1])=0=
\Hom_{\ltt{T}}(U_{k},\Theta(i+k+s)[1])= 0,$ it follows that
$\Hom_{\ltt{T}}(U,\Theta(i+k+s)[1])=0$ for any $s\in [-k,0].$ Thus $\Hom_{\ltt{T}}(U,\Theta(j)[1])=0$ for any $j\in[i,i+k+1]$. On the other hand, by \ref{filtradocero} (a), we have that
$\Hom_{\ltt{T}}(\Theta(i+k+1)^{a},U_{k})=0$ since $U_{k}\in
\F(\{\Theta(j)\;|\; j\in[i,i+k] \})$; also by
\ref{filtradocero} (c), we get that
$\Hom_{\ltt{T}}(\Theta(i+k+1)^{a},U_{k}[-1])=0$. Thus, by
\ref{triangulochido}, there exists a distinguished triangle
$$\eta':\quad \xymatrix{\Theta(i+k+1)^{d}\ar[r] & U_{k+1}\ar[r] &
U_{k}\ar[r] & \Theta(i+k+1)^{d}[1]}$$ with $d\leq a$ and
$U_{k+1}$ an indecomposable direct summand of $U$. By base change (see \ref{basecobase}), 
we have the following commutative diagram
$$\xymatrix{&  \Theta(i)[-1]\ar@{=}[r]\ar[d] & \Theta(i)[-1]\ar[d]\\
\Theta(i+k+1)^{d}\ar[r]\ar@{=}[d] &
V_{k+1}\ar[r]\ar[d]^{\mu_{k+1}} & V_{k}\ar[r]\ar[d]^{\mu_{k}} &
\Theta(i+k+1)^{d}[1]\ar@{=}[d]\\
\Theta(i+k+1)^{d}\ar[r] & U_{k+1}\ar[r]\ar[d] & U_{k}\ar[r]\ar[d] & \Theta(i+k+1)^{d}[1] \\
& \Theta(i)\ar@{=}[r] & \Theta(i),}$$
where the rows and columns are distinguished triangles. Using the fact that $V_{k} \in
\F(\{\Theta(j)\;|\;i<j \leq i+k\}),$ it follows by
\ref{extensionesclosed} that $V_{k+1} \in \F(\{\Theta(j)\;|\; i<j \leq i+k+1\})$. Moreover
$\Hom_{\ltt{T}}(U_{k+1},\Theta(j)[1])=0$ for $j\in[i,i+k+1],$
since $U_{k+1}$ is an indecomposable direct summand of $U$. Hence, the desired triangle is 
the first column of the preceding diagram, that is, $\xi_{k+1}$ is the triangle
$$\xymatrix{V_{k+1}\ar[r]^{\mu_{k+1}} &
U_{k+1}\ar[r] & \Theta(i)\ar[r] & V_{k+1}[1]}.$$
\end{dem}

\begin{teo}\label{sisexproy} Let $(\Theta,\leq )$ be a $\Theta$-system of
size $t,$ in an artin  triangulated $R$-category $\ltt{T}\!.$ Then,
there exists a unique, up to isomorphism, family $\Qb$ of objects in $\T$  such that
$(\Theta,\Qb,\leq )$ is a $\Theta$-projective system of size $t$ in
$\ltt{T}$.
\end{teo}
\begin{dem} Without lost of generality, we can assume that $\leq$ is the natural order on the set
$[1,t]$. For each $i<t$, we set $\eta_{i}:=\xi_{t-i}$ where
$\xi_{t-i}$ is the distinguished triangle of \ref{existriangulos}
$$\xi_{t-i}\;:\; \xymatrix{V_{t-i}\ar[r] & U_{t-i}\ar[r] &
\Theta(i)\ar[r] & V_{t-i}[1].}$$ Let $K(i):=V_{t-i}$ and
$Q(i):=U_{t-i}.$ Then, we have that $K(i)\in\F(\{\Theta(j)\;|\;
j>i\})$ and $\Hom_{\ltt{T}}(Q(i),\Theta(j)[1])=0$ for $j\geq
i$. From the triangle $\xi_{t-i}$, it follows that $Q(i)\in
\F(\{\Theta(j)\;|\; j\geq i\})$. By \ref{filtradocero} (b) and
(c), we conclude that $\Hom_{\ltt{T}}(Q(i),\Theta(r)[1])=0$ for
$r\leq i$ and $\Hom_{\ltt{T}}(Q(i),\Theta(r)[-1])=0$
$\forall\; r$. Therefore $Q(i)\in\,{}^{\perp}\Theta[-1]\cap{}^{\perp}\Theta[1]$ $\;\forall\,i.$ For $i=t$, we take the triangle $\eta_t$ as follows $$\xymatrix{0\ar[r] & \Theta(t)\ar[r]^{1}\ar[r] & \Theta(t)\ar[r] &
0}$$ and we set $Q(t):=\Theta(t)$  and $K(t):=0$, so this triangle
has the desired conditions. Finally, if there is another family $\Qb'$ such that
$(\Theta,\Qb',\leq )$ is a $\Theta$-projective system of size $t,$ then by \ref{Basico1Proj} 
we get that  $\Qb\simeq\Qb'.$ 
\end{dem}

\begin{lem}\label{exacto} Let $(\Theta,\Qb,\leq )$ be a $\Theta$-projective system of size $t,$ 
in a  triangulated $R$-category $\ltt{T}$. Then, the $R$-functor 
$\Hom_{\ltt{T}}(Q',-):\F(\Theta) \to \modu\,(R)$ is exact, for any $Q'\in \add\,(Q).$
\end{lem}
\begin{dem} Let  $\eta\;:\;A\to B\to C\to A[1]$ be  a distinguished triangle in $\F(\Theta)$ and 
$Q'\in \add\,(Q).$ Applying $\Hom_{\ltt{T}}(Q', -)$ to $\eta$, we get the exact
sequence
$$(Q',C[-1])\to (Q',A)\to (Q',B)\to (Q',C)\to (Q',A[1]).$$
Since $\Qb\subseteq{}^{\perp}\Theta[-1]\cap{}^{\perp}\Theta[1],$ it follows from 
\ref{filtradogeneral} that $\Hom_\T(Q',C[-1])=\Hom_\T(Q',A[1])=0.$ Thus, such a functor is exact.
\end{dem}

\begin{pro}\label{multindepen} Let $(\Theta,\Qb,\leq )$ be a $\Theta$-projective system of size $t,$ in a Hom-finite triangulated $R$-category $\ltt{T}$. Then, the following statements hold.
\begin{enumerate}
\item [(a)] For any $M\in \F(\Theta)$ the filtration multiplicity $[M:\Theta(i)]_\xi$ of $\Theta(i)$ in $M$ doest not depend on the given $\Theta$-filtration
$\xi$ of $M$ and hence it will be denoted by $[M:\Theta(i)]$. In particular $\ell_{\Theta}(M)=\sum_{i=1}^{t}[M:\Theta(i)]$.
\item [(b)] $Q(i)\not\simeq Q(j)$ if $i\neq j.$
\end{enumerate}
\end{pro}
\begin{dem} (a) Consider a  $\Theta$-filtration $\xi$ of $M\in \F(\Theta)$
$$\xi=\{\xi_{l}\;:\; \xymatrix{M_{l-1}\ar[r] & M_{l}\ar[r] &
\Theta(j_{l})\ar[r] & M_{l-1}[1]}\}_{l=0}^{n},$$ where 
$M_{-1}=0=\Theta(j_{0})$, $j_{l}\in [1,t]$ for $l\geq 1,$ and
$M_{n}=M$. Applying the functor $\Hom_{\ltt{T}}(Q(i),-)$ to each
triangle $\xi_{j},$ and by setting $\langle X,Y\rangle:=\ell_{R}(\Hom_{\ltt{T}}(X,Y)),$ we get 
the following equalities
\begin{align*} \langle Q(i),M_{1}\rangle & =\langle
Q(i),0\rangle + \langle Q(i),\Theta(j_{1})\rangle,\\
\langle Q(i),M_{2}\rangle & =\langle
Q(i),\Theta(j_{1})\rangle + \langle Q(i),\Theta(j_{2})\rangle,\\
\langle Q(i),M_{3}\rangle  & =\langle Q(i),M_{2}\rangle +
\langle Q(i),\Theta(j_{3})\rangle,\\
 & \vdots\\
\langle Q(i),M\rangle & =\langle Q(i),M_{n-1}\rangle + \langle
Q(i),\Theta(j_{n})\rangle.
\end{align*}
Let $c_i:=\langle Q(i),M\rangle
=\sum_{j=1}^{t}[M:\Theta(j)]_{\xi}\langle Q(i),\Theta(j)\rangle.$
Consider the matrix $D:=(d_{ij}),$ where $d_{ij}:=\langle Q(i),\Theta(j)\rangle$. By
\ref{anulasuperior} (b), we have that $D$ is an upper triangular matrix with $d_{ii}\neq
0$ $\;\forall\,i,$ and thus $\text{det}\,(D)\neq 0$. By using the column vectors
$X:=([M:\Theta(1)]_{\xi},[M:\Theta(2)]_{\xi}
,\cdots,[M:\Theta(t)]_{\xi})^{t}$ and $C:=(c_1,c_2,\cdots,c_t)^{t},$ the above equalities can 
be written as a matrix equation $D\cdot X=C$. Since $\text{det}\,(D)\neq 0$, we obtain that 
$X=D^{-1}\cdot C,$ and hence $[M:\Theta(j)]_{\xi}$ only depends on the numbers 
$c_i=\langle Q(i),M\rangle$ and $d_{ij}=\langle Q(i),\Theta(j)\rangle.$ 
\

(b) Let $i\neq j.$ We can assume that $j>i.$ Then, by (a) and \ref{epss} (PS5), it follows that 
$[Q(i):\Theta(i)]=1$ and $[Q(j):\Theta(i)]=0,$ and thus $Q(i)\not\simeq Q(j).$
\end{dem}

\begin{defi} Let $(\Theta,\Qb,\leq )$ be a $\Theta$-projective system of size $t,$ in a 
Hom-finite triangulated $R$-category $\ltt{T}$. The {\bf $\Theta$-support} of $M\in\F(\Theta),$ is 
the set $$\Supp_\Theta(M):=\{i\in[1,t]\;|\; [M:\Theta(i)]\neq 0\}.$$
\end{defi}

For $0\neq M\in\F(\Theta),$ let $\maxi\,(M)$ denote the maximum of $\Supp_\Theta(M)$ with respect to the linear order $\leq,$ and similarly,  $\mini\,(M)$ denote the minimum of $\Supp_\Theta(M)$ with respect to the linear order $\leq.$ Finally, we set $\maxi\,(0):=-\infty$ and $\mini\,(0):=+\infty.$

\begin{teo}\label{proyecrelativo} Let $(\Theta,\Qb, \leq)$ be a $\Theta$-projective system 
of size $t$, in a Hom-finite triangulated $R$-category $\ltt{T},$ and let $M\in \F(\Theta)$ and 
$i:=\mini\,(M)$. Then, there exists a distinguished triangle in $\T$
$$\xymatrix{N\ar[r] & Q_{0}(M)\ar[r]^{\varepsilon_{M}} &
M\ar[r] & N[1]}$$ satisfying the following conditions:
\begin{enumerate}
\item [(a)] $N\in \F(\Theta)$ and $Q_{0}(M)\in\add\,(\bigoplus_{j\geq i}\,Q(j)),$
\item [(b)] $\mini\,(M)<\mini\,(N)$ if $M\neq 0,$ 
\item [(c)] $\varepsilon_{M}:Q_{0}(M)\to M$ is a $\Pro(\Theta)$-precover of $M.$
\end{enumerate}
\end{teo}
\begin{dem} If $M=0,$ the zero distinguished triangle $0\to 0\to 0\to 0$ is the desired one. Without lost of generality, it can be assumed that 
$\leq$ is the natural order on the set $[1,t].$
\

 Let $M\neq 0.$ Then, by \ref{anulasuperior} (a) and \ref{filtracionordenada} (c), there is a distinguished triangle
$$\xymatrix{N\ar[r]^{\varphi} & M\ar[r]^{\psi} &
\Theta(i)^{m_{i}}\ar[r] & N[1]}$$ with $N\in \F(\Theta)$
and $\mini\,(M)<\mini\,(N)$. We proceed by reverse induction on $i=\mini\,(M).$ If $i=\mini\,(M)=t$, we have that $N=0$ and hence the desired triangle is $0\to \Theta(t)^{m_{i}}\to 
\Theta(t)^{m_{i}}\to  0,$  since $Q(t)\simeq \Theta(t)$.\\
Let $i=\mini\,(M)<t$. If $N=0$, we have that $M=\Theta(i)^{m_{i}}$ and thus the following distinguished triangle (see \ref{epss} (PS5)) is the 
desired one
$$\xymatrix{K(i)^{m_{i}}\ar[r] & Q(i)^{m_{i}}\ar[r]^{\beta_i^{m_i}} & \Theta(i)^{m_{i}}\ar[r] &
 K(i)^{m_{i}}[1].}$$ Suppose that $N\neq 0.$ Since $i=\mini\,(M)<\mini\,(N)$, by induction, there is a
distinguished triangle
$$\xymatrix{N'\ar[r] & Q_{0}(N)\ar[r]^{\varepsilon_{N}} & N\ar[r] & N'[1]}$$
such that $i<\mini\,(N)<\mini\,(N')=:i'$, $Q_{0}(N)\in\text{add}\,(\bigoplus_{j\geq i'}Q(j))$ and $\varepsilon_N:Q_{0}(N) \to
N$ is a $\Pro(\Theta)$-precover of $N$. By base change (see \ref{basecobase}), we obtain the following commutative diagram in $\T$
$$\xymatrix{& K(i)^{m_{i}}\ar@{=}[r]\ar[d] &  K(i)^{m_{i}}\ar[d]\\
N\ar[r]^{i_1}\ar@{=}[d] & E\ar[r]^{p_2}\ar[d]^{\theta} &
Q(i)^{m_{i}}\ar[r]\ar[d]^{\beta_{i}^{m_i}}
& N[1]\ar@{=}[d]\\
N\ar[r]^{\varphi} & M\ar[r]^{\psi}\ar[d] &
\Theta(i)^{m_{i}}\ar[r]\ar[d] & N[1]\\
&  K(i)^{m_{i}}[1]\ar@{=}[r] & K(i)^{m_{i}}[1].}$$ Since
$N\in \F(\Theta)$, we have that  $\Hom_{\ltt{T}}(Q(i),N[1])=0$.
Thus, the first row, in the diagram above, splits. So, there
is $i_{2}:Q(i)^{m_{i}}\to E$ such that $\beta_i^{m_i}=\psi\theta
i_{2}$. Define $\alpha:=\theta i_{2}$  and
$\varepsilon:=(\varphi\varepsilon_{N},\alpha):Q_{0}(N)\bigoplus
Q(i)^{m_{i}} \fl M$. Hence, we get the following commutative diagram
$$\xymatrix{Q_{0}(N)\ar[r]^{j_1}\ar[d]^{\varepsilon_{N}}\ar@{}[dr] &
Q_{0}(M)\ar[r]^{\pi_2}\ar[d]^{\varepsilon}\ar@{}[dr] &
Q(i)^{m_{i}}\ar[d]^{\beta_{i}^{m_{i}}}\ar[r]^{0}\ar@{}[dr]
& Q_{0}(N)[1]\ar[d]^{\varepsilon_{N}[1]}\\
N\ar[r]^{\varphi} & M\ar[r]^{\psi} & \Theta(i)^{m_{i}}\ar[r] &
N[1]}$$ with $Q_{0}(M):=Q_{0}(N)\bigoplus Q(i)^{m_{i}},$ where the rows are 
distinguished triangles, $j_{2}:=\left(\begin{array}{c}
1\\
0
\end{array}\right)$ and $\pi_{2}:=\left(\begin{array}{cc}
0 & 1
\end{array}\right)$. Since $\Qb\subseteq {}^{\perp}\Theta[-1]$ and 
$Q_{0}(N)\in \add\,(\bigoplus_{j\geq i'}Q(j)),$  we 
conclude that $\Hom_{\ltt{T}}(Q_{0}(N),\Theta(i)^{m_{i}}[-1])=0.$ Thus, 
by \ref{SnakeLema}, we obtain the following diagram in $\T$ 
$$\xymatrix{N'\ar[r]\ar[d] &
P[-1]\ar[r]\ar[d] & K(i)^{m_{i}}\ar[r]\ar[d] &
 N'[1]\ar[d]\\
Q_{0}(N)\ar[r]^{j_{1}}\ar[d]^{\varepsilon_{N}} &
Q_{0}(M)\ar[r]^{\pi_{2}}\ar[d]^{\varepsilon} &
Q(i)^{m_{i}}\ar[r]\ar[d]^{\beta_{i}^{m_{i}}} & 
Q_{0}(N)[1]\ar[d]^{\varepsilon_{N}[1]}\\
N\ar[r]^{\varphi}\ar[d] & M\ar[r]^{\psi}\ar[d]^{\gamma} & \Theta(i)^{m_{i}}\ar[r]\ar[d]\ar@{}[dr]|{IX} &  N[1]\ar[d]\\
 N'[1]\ar[r] & P\ar[r]  & K(i)^{m_{i}}[1]\ar[r] &
N'[2].}$$ where the rows and columns, in the diagram above, are distinguished triangles and all squares commute, except the one marked with $IX,$ which anti-commutes. We claim that the
following distinguished triangle in $\T$
$$\xymatrix{P[-1]\ar[r] & Q_{0}(M)\ar[r]^{\varepsilon} & M\ar[r] & P}$$
is the desired one. Indeed, we have that $Q_{0}(M) \in \text{add}\,(\bigoplus_{j\geq i}Q(j))$ since
$Q_{0}(N)\in \text{add}\,(\bigoplus_{j\geq i'}Q(j))$ with
$i<\mini\,(N)<\mini\,(N')=i'$ and $Q(i)^{m_{i}}\in
\F(\{\Theta(j)\;|\; j\geq i\})$. By considering the first row from
the last diagram, it follows that $P[-1]\in\F(\{\Theta(j)\;|\; j>i\})$ since $K(i)^{m_{i}}\in \F(\{\Theta(j)\;|\; j>i\})$ and
$N'\in \F(\{\Theta(j)\;|\; j>i'\})$ with $i<i'$. Therefore
$i=\mini\,(M)<\mini\,(P[-1])$.\\
Finally, we show that $\varepsilon$ is a $\Pro(\Theta)$-precover of $M$. Indeed, 
let $h:X \to M$ be a morphism in $\T$ with $X\in\Pro(\Theta),$ and consider the
morphism $\gamma h:X \to P.$ Since $P[-1]\in\F(\Theta)$, we have that 
$\gamma h=0$ and so there is  $h':X \to Q_{0}(M)$  such
that $h=\varepsilon h'$. Then, the morphism $\varepsilon$ is a $\Pro(\Theta)$-precover of $M$, 
proving the result. 
\end{dem}

\begin{cor}\label{addQ} Let $(\Theta,\Qb, \leq)$ be a $\Theta$-projective system 
of size $t$, in a Hom-finite triangulated $R$-category $\ltt{T}.$ Then 
$$\add\,(Q)=\F(\Theta)\cap \Pro(\Theta).$$
\end{cor}
\begin{dem} It is clear that $\add\,(Q)\subseteq\F(\Theta)\cap \Pro(\Theta).$ Let $M\in
\F(\Theta)\cap \Pro(\Theta).$ Then, by \ref{proyecrelativo}, there is a distinguished triangle 
in $\T$
$$\eta\;:\; \xymatrix{N\ar[r] & Q_{0}(M)\ar[r]^{\varepsilon_{M}} &
M\ar[r] & N[1],}$$
where $Q_{0}(M)\in\add\,(Q)$ and $N\in\F(\Theta).$ Thus, the triangle $\eta$ splits and then 
$M\in\add\,(Q);$ proving the result.
\end{dem} 

\section{The standardly stratified algebra associated to a $\Theta$-projective system}

\begin{teo}\label{TeoEquFilt} Let $(\Theta,\Qb,\leq)$
be a $\Theta$-projective system of size $t,$ in an artin triangulated $R$-category $\ltt{T},$ and 
let $A:=\End_\T(Q)^{op},$ $e_Q:=\Hom_\T(Q,-):\T\to\modu\,(A)$ and ${}_AP(i):=e_Q(Q(i))$ for each $i\in[1,t].$ 
Then, the following statements holds.
\begin{enumerate}
\item[(a)] The family ${}_AP:=\{{}_AP(i)\;|\;i\in[1,t]\}$ is a representative set of the indecomposable projective $A$-modules. In particular, $A$ is basic and $rk\,K_0(A)=t.$
\item[(b)] $e_Q(\Theta(i))\simeq {}_A\Delta(i)$ $\;\forall\,i\in[1,t],$  where ${}_A\Delta$ is 
computed by using ${}_AP$ and the given order $\leq$ on $[1,t].$
\item[(c)] $(A,\leq)$ is a standardly stratified algebra, that is, $\proj\,(A)\subseteq 
\F({}_A\Delta).$
\item[(d)] The restriction $e_Q:\F(\Theta)\to \F({}_A\Delta)$ is an exact equivalence of 
$R$-categories.
\end{enumerate}
\end{teo}
\begin{dem} By \ref{exacto}, we know that $e_Q=\Hom_\T(Q,-)|_{\F(\Theta)}:\F(\Theta)\to\modu\,(A)$ 
is an exact functor.
\

 (a) It follows by \ref{ProAusl} and \ref{multindepen} (b).
\

(b) and (c) Let $i\in[1,t].$ By \ref{epss} (PS5) and \ref{proyecrelativo}, we have two 
distinguished triangles in $\T$
$$\begin{matrix}\xymatrix{\eta_i\;:\;K(i)\ar[r]^{\alpha_i} & Q(i)\ar[r] & \Theta(i)\ar[r] & K(i)[1],}\\{}\\
\eta'_i\;:\; \xymatrix{K'\ar[r] & Q'\ar[r]^{\lambda_i} & K(i)\ar[r] & K'[1],}
   \end{matrix}$$
where $K(i), K'\in\F(\{\Theta(j)\;|\;j>i\}$ and $Q'\in\add\,(\oplus_{j>i}\;Q(j)).$ Applying the functor $e_Q=\Hom_\T(Q,-)$ to the triangles $\eta_i$ and $\eta'_i,$ we get the following exact 
sequence in $\modu\,(A)$
$$\varepsilon_i\;:\;\xymatrix{e_Q(Q')\ar[r]^{e_Q(\gamma_i)} & {}_AP(i)\ar[r] & e_Q(\Theta(i))\ar[r] & 0,}$$
where $\gamma_i:=\alpha_i\lambda_i.$ We assert that $$\Ima\,(e_Q(\gamma_i))=\Tr_{\oplus_{j>i}\;{}_AP(j)}\,({}_AP(i)).$$
Indeed, using the fact that $e_Q(Q')\in\add\,(\oplus_{j>i}\;{}_AP(j)),$ it follows 
that\\ $\Ima\,(e_Q(\gamma_i))\subseteq \Tr_{\oplus_{j>i}\;{}_AP(j)}\,({}_AP(i)).$ In order to see the other inclusion, let $j>i$ and consider a morphism $f:{}_AP(j)\to {}_AP(i).$ By \ref{epss} (PS3) and \ref{ProAusl} (c), 
we conclude that $\Hom_A({}_AP(j),e_Q(\Theta(i))=0$ and hence $f$ factorizes trough $e_Q(\gamma_i);$ proving our assertion. Finally, by this assertion and the exact sequence 
$\varepsilon_i,$ we obtain (b) and (c).
\

(d) Since $e_Q:\F(\Theta)\to\modu\,(A)$ is an exact functor, it remains to prove that 
$e_Q:\F(\Theta)\to \F({}_A\Delta)$ is full, faithful and a dense functor.
\

Let $M\in \F(\Theta).$ We prove, by induction on $\ell_{\Theta}(M),$ that $e_{Q}(M)\in
\F({}_{A}\Delta)$. If $\ell_{\Theta}(M)\leq 1$, then $M=\Theta(i)^{m_i}$ for some $i,$ and hence 
by (a) it follows that $e_{Q}(M)\in\F({}_{A}\Delta)$.\\
Let $\ell_{\Theta}(M)>1.$ Then, from \ref{anulasuperior} (a) and \ref{filtracionordenada} (c), 
there is a distinguished triangle $N\to M\to \Theta(i)^m\to N[1]$ in $\F(\Theta)$ such that 
$\ell_{\Theta}(N)<\ell_{\Theta}(M).$ Therefore, by induction and since $\F({}_{A}\Delta)$
is closed under extensions, we conclude that $e_{Q}(M)$ belongs to $\F({}_{A}\Delta);$ proving 
that $\Ima\,(e_{Q}|_{\F(\Theta)})\subseteq\F({}_{A}\Delta).$
\

Now, we prove that $e_Q:\F(\Theta)\to \F({}_A\Delta)$ is full and faithful. Indeed, let $M,N\in
\F(\Theta).$ By \ref{proyecrelativo} and the exactness of the functor $e_Q,$ we get an exact 
sequence $\varepsilon\;:\; e_Q(Q_1)\to e_Q(Q_0)\to M\to 0$ in $\modu\,(A)$ such that $Q_0,Q_1\in 
\add\,Q.$ From $\varepsilon,$ we obtain the following exact and
commutative diagram 
$$\xymatrix{0\ar[r] & _{\ltt{T}}(M,N)\ar[r]\ar[d]^{\alpha_1} &
_{\ltt{T}}(Q_{0},N)\ar[r]\ar[d]^{\alpha_2} &
_{\ltt{T}}(Q_{1},N)\ar[d]^{\alpha_3}\\
0\ar[r] & _{A}(e_{Q}(M),e_{Q}(N))\ar[r] &
_{A}(e_{Q}(Q_{0}),e_{Q}(N))\ar[r] &
_{A}(e_{Q}(Q_{1}),e_{Q}(N))}$$ where $\alpha_{2}$ and
$\alpha_{3}$ are isomorphism (see \ref{ProAusl} (c)). Thus, by using the so-called Five's Lemma, it follows 
that $\alpha_{1}$ is an isomorphism; proving that $e_{Q}$ is full and faithful.
\

Finally, we see that $e_Q:\F(\Theta)\to \F({}_A\Delta)$ is dense.  Indeed, let 
$M\in\F({}_A\Delta).$ We proceed by induction on the ${}_{A}\Delta$-length  
$\ell_{{}_{A}\Delta}(M).$ If $\ell_{{}_{A}\Delta}(M)=1$ then 
$M\simeq{}_{A}\Delta(i)\simeq e_{Q}(\Theta(i))$ for some $i.$\\
Let $\ell_{{}_{A}\Delta}(M)>1.$ Then, there is an exact sequence in $\modu\,(A)$
$$\xymatrix{0\ar[r] & {}_{A}\Delta(i)\ar[r] & M\ar[r] &
M/{}_{A}\Delta(i)\ar[r] & 0,}$$ where
$\ell_{{}_{A}\Delta}(M/{}_{A}\Delta(i))=\ell_{{}_A\Delta}(M)-1$ for some $i.$ So, by
induction, there exists $Z\in \F(\Theta)$ such that
$e_{Q}(Z)\simeq M/{}_{A}\Delta(i)$. Moreover, by \ref{proyecrelativo}, there
is a distinguished triangle $\eta_{Z}\;:\;Z'\stackrel{u}{\to} Q_{0}(Z)\stackrel{\varepsilon_{Z}}
{\to} Z\to Z'[1]$ in $\F(\Theta),$ with $Q_{0}(Z)\in\add\,(Q);$ and thus, we get the following 
exact and commutative diagram in $\modu\,(A)$
$$\xymatrix{ & & 0\ar[d] & 0\ar[d]\\
& & e_{Q}(Z')\ar@{=}[r]\ar[d]^{\mu} & e_{Q}(Z')\ar[d]^{e_{Q}(u)}\\
\eta\;:\; 0\ar[r] & {}_{A}\Delta(i)\ar[r]^{i_{1}}\ar@{=}[d] &
C\ar[r]^{p_2}\ar[d]^{\lambda} &
e_{Q}(Q_{0}(Z))\ar[r]\ar[d]^{e_{Q}(\varepsilon_{Z})} & 0\\
0\ar[r] & {}_{A}\Delta(i)\ar[r] & M\ar[r]\ar[d] &
e_{Q}(Z)\ar[r]\ar[d] & 0\\
& & 0 & 0.}$$ 
Since $e_{Q}(Q_{0}(Z))\in\proj\,(A),$ the exact sequence $\eta$ splits and hence $C=\,
{}_{A}\Delta(i)\bigoplus e_{Q}(Q_{0}(Z))\simeq
e_{Q}(\Theta(i)\bigoplus Q_{0}(Z))$,
$i_{1}=\left(\begin{array}{c}
1\\
0
\end{array}\right)$ and $p_{2}=(0,1)$. That is $\mu=\left(\begin{array}{c}
\varphi\\
e_{Q}(u) \end{array}\right)$ with $\varphi:e_{Q}(Z') \to e_{Q}(\Theta(i))$. Since the restriction 
$e_{Q}\mid_{\F(\Theta)}$ is full, there exists  $h:Z' \to \Theta(i)$ such that 
$e_{Q}(h)=\varphi$ and hence $\mu=e_{Q}(\psi),$ where $\psi:=\left(\begin{array}{c}
h\\
u
\end{array}\right)$. Then, by completing $\psi$ to a distinguished triangle and from \ref{SnakeLema}, we get the following commutative diagram 
$$\xymatrix{& \Theta(i)\ar@{=}[r]\ar[d] &  \Theta(i)\ar[d]\\
Z'\ar[r]^(.4){\psi}\ar@{=}[d] & \Theta(i)\bigoplus
Q_{0}(Z)\ar[r]\ar[d]^{\pi_2} & X\ar[r]\ar[d]^{\alpha}
& Z'[1]\ar@{=}[d]\\
Z'\ar[r]^{u} & Q_{0}(Z)\ar[r]^{\varepsilon_{Z}}\ar[d] &
Z\ar[r]\ar[d] & Z'[1]\\
&  \Theta(i)[1]\ar@{=}[r] & \Theta(i)[1],}$$
where the rows and columns are distinguished triangles and $\pi_2:=(0,1)$. Observe that 
$X\in \F(\Theta)$ since $\F(\Theta)$ is closed under
extensions. Thus, by applying $e_{Q}$ to the first row, in the diagram above, we get 
the exact sequence
$$\xymatrix{0\ar[r] & e_{Q}(Z')\ar[r]^(.4){e_{Q}(\psi)} & e_{Q}(\Theta(i)\bigoplus
Q_{0}(Z))\ar[r] & e_{Q}(X)\ar[r] & 0.}$$ But $e_{Q}(\psi)=\mu$ and hence 
$e_{Q}(X)\simeq \Coker\,(\mu)=M;$ proving that $e_Q:\F(\Theta)\to \F({}_A\Delta)$ is dense.
\end{dem}

\begin{pro}\label{ProjCover} Let $(\Theta,\Qb,\leq)$ be a $\Theta$-projective system of 
size $t,$ in an artin triangulated $R$-category $\ltt{T}.$ Then, the following statements hold.
\begin{itemize}
\item[(a)] $\F(\Theta)$ is closed under extensions and direct summands.
\item[(b)] $\Theta(i)$ is indecomposable for each $i\in[1,t].$
\item[(c)] For any object $M\in \F(\Theta),$ there exists a distinguished triangle 
$Z\to Q_M\to M\to Z[1]$ in $\F(\Theta)$ such that $Q_M\to M$ is an $\add\,(Q)$-cover 
of $M,$ and $\mini\,(M)<\mini\,(Z)$ if $M\neq 0.$
\end{itemize}  
\end{pro}
\begin{dem} Let $A:=\End_\T(Q)^{op}.$ We know by \ref{TeoEquFilt} that: 
$e_Q:\F(\Theta)\to \F({}_A\Delta)$ is an exact equivalence, $(A,\leq)$ is an standardly stratified algebra and $e_Q(\Theta(i))\simeq {}_A\Delta(i)$ $\;\forall\,i.$ Since ${}_A\Delta(i)$ 
is indecomposable and $\End_\T(\Theta(i))\simeq \End_A({}_A\Delta(i)),$ it follows (b). To prove 
(a), we use the well-known fact that $\F({}_A\Delta)$ is closed under direct summands (see 
\cite{AHLU}). Indeed, this property can be carried back to $\F(\Theta)$ by using the equivalence 
$e_Q$ and that both $\T$ and $\modu\,(A)$ are Krull-Schmidt categories.
\

Finally,  since $\F({}_A\Delta)$ is a resolving subcategory of $\modu\,(A)$ (see \cite{AHLU}), by 
using \ref{proyecrelativo} and the exact equivalence  $e_Q:\F(\Theta)\to \F({}_A\Delta),$ 
we get (c).
\end{dem}

\begin{cor}\label{epss=ss} Let $(\Theta,\Qb,\leq)$ be a $\Theta$-projective system of 
size $t,$ in an artin triangulated $R$-category $\ltt{T},$ and let $\Kb:=\{K(i)\}_{i=1}^t$ where,
for each $i,$ $K(i)$ is the object appearing in \ref{epss} (PS5). If $\Hom_\T(\Kb[2],\Theta)=0$ 
then $(\Theta,\leq)$ is a $\Theta$-system of size $t,$ in $\ltt{T}.$
\end{cor}
\begin{dem} It follows from \ref{ProjCover} (b) and \ref{anulasuperior} (a), (c).
\end{dem}

\begin{cor}\label{PropThetaSys}  Let $(\Theta,\leq)$ be a $\Theta$-system of 
size $t,$ in an artin triangulated $R$-category $\ltt{T}.$ Then, the following statements hold.
\begin{itemize}
\item[(a)] $\F(\Theta)$ is closed under extensions and direct summands.
\item[(b)] There is a unique, up to isomorphism, $\Theta$-projective system $(\Theta,\Qb,\leq)$ of size $t,$ which is associated to 
the $\Theta$-system $(\Theta,\leq).$
\item[(c)] For any object $M\in \F(\Theta),$ there exists a distinguished triangle 
$Z\to Q_M\to M\to Z[1]$ in $\F(\Theta)$ such that $Q_M\to M$ is an $\add\,(Q)$-cover 
of $M,$ and $\mini\,(M)<\mini\,(Z)$ if $M\neq 0.$
\item[(d)] $\F(\Theta)\cap\Pro(\Theta)=\add\,(Q).$
\end{itemize}
\end{cor} 
\begin{dem} It follows from \ref{sisexproy}, \ref{ProjCover} and \ref{addQ}.
\end{dem}
\vspace{.2cm}

The previous results can be seen also under the light of the so-called cotorsion pairs in the sense of Iyama-Nakaoka-Yoshino
(see \cite{IY} and \cite{N}). Such cotorsion pairs are studied extensively in relation with cluster tilting categories, $t$-structures 
and co-$t$-structures.

\begin{defi} A pair $(\X,\Y)$ of subcategories in a triangulated category $\T$ is called a cotorsion pair if the following conditions 
hold.
\begin{itemize}
 \item[(a)] $\X$ and $\Y$ are closed under direct summands in $\T.$
 \item[(b)] $\Hom_\T(\X,\Y)=0\:$ and $\;\T=\X*\Y[1].$
\end{itemize}
The core of the cotorsion pair $(\X,\Y)$ is the subcategory $\X\cap\Y.$
\end{defi}

\begin{cor} \label{CotorsionThetaSys1}  Let $(\Theta,\leq)$ be a $\Theta$-system of 
size $t,$ in an artin triangulated $R$-category $\ltt{T}.$ Then, the pairs $(\Pro(\Theta),\F(\Theta))$ and $(\F(\Theta),\I(\Theta))$ are 
cotorsion pairs. 
\end{cor}
\begin{dem} Since $\Pro(\Theta):={}^\perp\F(\Theta)[1]$ and $\I(\Theta):=\F(\Theta)^\perp[-1],$ it follows that these classes are 
closed under direct summands in $\T.$ Furthermore, by \ref{PropThetaSys} (a), we also know that $\F(\Theta)$ is closed 
under direct summands in $\T.$ Finally, from \ref{TeoRingel2}, we get that $\Pro(\Theta)*\F(\Theta)[1]=\T=\F(\Theta)*\I(\Theta)[1].$ 
\end{dem}

\begin{rk}\label{CotorsionThetaSys2} Let $(\Theta,\leq)$ be a $\Theta$-system of size $t,$ in an artin triangulated $R$-category $\ltt{T}.$ Observe that, by \ref{InjVsProj}, \ref{TeoEquFilt} and \ref{PropThetaSys}, 
the cotorsion pairs  $(\Pro(\Theta),\F(\Theta))$ and $(\F(\Theta),\I(\Theta))$ have the following properties. Their cores are determined by, respectively, the $\Theta$-projective system $(\Theta,\Qb,\leq)$ and the $\Theta$-injective system $(\Theta,\Yb,\leq)$ as follows:
  $$\F(\Theta)\cap\Pro(\Theta)=\add\,(Q)\quad\text{and}\quad \F(\Theta)\cap\I(\Theta)=\add\,(Y).$$
Moreover, the endomorphism algebras $\End_\T(Q)^{op}$ and $\End_\T(Y)$ are standardly stratified algebras.
\end{rk}

\section{The bounded derived category $\D^b(\F(\Theta))$}

We recall that an exact category is  and additive category $\A$ endowed with a class $\E$ of pairs $M\stackrel{i}{\to} E
\stackrel{p}{\to}N$ in $\A$ closed under isomorphisms and satisfying a list of axioms \cite{Q, Keller2}. An exact category 
$(\A,\E)$ is saturated if every idempotent in $\A$ splits and so, in this case (see, for example in \cite{Keller3, Neeman}), there 
exists the bounded derived category $\D^b(\A).$ 
\

Let $(A,\leq)$ be an standardly stratified algebra. It is well-known that $\F({}_A\Delta)$ is an additive category, which is closed 
under extensions and every idempotent in $\F({}_A\Delta)$ splits. Consider the class $\Ex\,({}_A\Delta)$ of all pairs 
$M\stackrel{i}{\to} E \stackrel{p}{\to}N$ in $\F({}_A\Delta)$ such that $0\to M\stackrel{i}{\to} E \stackrel{p}{\to}N\to 0$ is 
an exact sequence in $\modu\,(A).$ Then, the pair $(\F({}_A\Delta),\Ex\,({}_A\Delta))$ is an exact category, and since 
it is saturated, there exists the bounded derived category $\D^b(\F({}_A\Delta)).$ We denote by $\D^b(A)$ to the bounded derived 
category of the abelian category $\modu\,(A).$

\begin{lem}\label{R1} Let $(A,\leq)$ be an standardly stratified algebra. Then $$\D^b(\F({}_A\Delta))\simeq\D^b(A)$$ as 
triangulated categories. 
\end{lem}
\begin{dem} Since $(A,\leq)$ is an standardly stratified algebra, it follows by \cite{AHLU} 
that $\F({}_A\Delta)$ is a resolving subcategory of $\modu\,(A).$ In particular, for any 
$M\in\F({}_A\Delta),$ there is an exact sequence $0\to M'\to P_0(M)\to M\to 0$ lying in 
$\F({}_A\Delta)$ and such that $P_0(M)\to M$ is the projective cover of $M.$ Hence, the 
construction outlined in \cite[Section 2]{BM} give us an equivalence 
$$R_{>-\infty}:\D^b(\F({}_A\Delta))\to K^{-,b}(\proj\,(A))$$ as triangulated categories, 
where $R_{>-\infty}:=\cdots R_{-2}R_{-1}R_0R_{>0}.$ So, 
the lemma follows, since $K^{-,b}(\proj\,(A))\simeq \D^b(A)$ as triangulated categories. 
\end{dem}

\begin{defi} Let $\Theta$ be a class of objects in a triangulated category $\T.$ We denote by 
$\Ex\,(\Theta)$ to the class of all the pairs $M\stackrel{i}{\to} E \stackrel{p}{\to}N$ in $\F(\Theta)$ 
admitting a morphism $q:N\to M[1]$ such that 
$M\stackrel{i}{\to} E \stackrel{p}{\to}N\stackrel{q}{\to}M[1]$ is a distinguished triangle in $\T.$
\end{defi}

\begin{teo}\label{R2} Let $(\Theta,\Qb,\leq)$ be a $\Theta$-projective system of size $t,$ in 
an artin triangulated $R$-category $\ltt{T}.$ Consider $A:=\End_\T(Q)^{op}$ and the 
functor $e_Q:=\Hom_\T(Q,-):\T\to \modu\,(A),$ 
where $Q:=\oplus_{i=1}^t\,Q(i).$ 
Then, the following statements hold.
\begin{enumerate}
\item[(a)] The pair $(\F(\Theta),\Ex\,(\Theta))$ is an exact and Krull-Schmidt category. Moreover, 
the equivalence $e_Q:\F(\Theta)\to\F({}_A\Delta)$ satisfies that 
$e_Q(\Ex\,(\Theta))=\Ex\,({}_A\Delta).$ 
\item[(b)] The derived functor $R\Hom_\T(Q,-):\D^b(\F(\Theta))\to \D^b(A)$ is an equivalence of triangulated categories.
\item[(c)] $\D^b(\F(\Theta))\simeq K^{-,b}(\add\,(Q))$ as triangulated categories. 
\end{enumerate}
\end{teo}
\begin{dem} (a) By \ref{extensionesclosed} and \ref{ProjCover} (a), we know that $\F(\Theta)$ is 
closed under extensions and direct summands in the artin triangulated $R$-category $\T.$ Thus 
$\F(\Theta)$ is an additive and Krull-Schmidt category. Consider the class $\Ex_Q\,(\Theta)$ of 
all the pairs $M\stackrel{i}{\to} E \stackrel{p}{\to}N$ in $\F(\Theta)$ satisfying that 
$e_Q(M)\stackrel{e_Q(i)}{\to} e_Q(E) \stackrel{e_Q(p)}{\to}e_Q(N)$ belongs to $\Ex\,({}_A\Delta).$ 
Since $e_Q:\F(\Theta)\to\F({}_A\Delta)$ is an exact equivalence of $R$-categories 
(see \ref{TeoEquFilt}) and $\F({}_A\Delta)$ is closed under extensions, it follows that 
$\Ex\,(\Theta)\subseteq\Ex_Q\,(\Theta)$ and also that the pair $(\F(\Theta),\Ex_Q\,(\Theta))$ is 
an exact category. It remains to see that $\Ex_Q\,(\Theta)\subseteq \Ex\,(\Theta).$
\

Let $M\stackrel{i}{\to} E \stackrel{p}{\to}N$ be in $\Ex_Q\,(\Theta).$ Consider a distinguished 
triangle of the form $\eta\;:\;M\stackrel{i}{\to} E \stackrel{\lambda}{\to}C\stackrel{\omega}
{\to}M[1].$ We assert that $e_Q(\lambda):e_Q(E)\to e_Q(C)$ is surjective. In order to see that, we 
apply the functor $e_Q$ to the triangle $\eta;$ and then we get the exact sequence 
$e_Q(E)\stackrel{e_Q(\lambda)}{\to} e_Q(C) \stackrel{e_Q(\omega)}{\to}e_Q(M[1]).$ But $e_Q(M[1])=0$ 
(see \ref{epss} (PS4)) and so  $e_Q(\lambda):e_Q(E)\to e_Q(C)$ is surjective. On the other hand, since 
$e_Q(p)e_Q(i)=0,$ it follows that $pi=0$ and hence there exists $p':C\to N$ such that $p'\lambda =p.$ Therefore, 
we get the following exact en commutative diagram in $\modu\,(A)$
\[\xymatrix{0\ar[r] & e_Q(M)\ar[r]^{e_Q(i)}\ar@{=}[d] &  e_Q(E)\ar[r]^{e_Q(\lambda)}\ar@{=}[d] & e_Q(C)\ar[r]\ar[d]^{e_Q(p')} & 0\\
   0\ar[r] & e_Q(M)\ar[r]^{e_Q(i)} &  e_Q(E)\ar[r]^{e_Q(p)} & e_Q(N)\ar[r] & 0.           }\]
Hence $e_Q(p')$ is an isomorphism in $\modu\,(A)$ and then $p':C\to N$ is an isomorphism in $\T.$ So the triangle 
$M\stackrel{i}{\to} E \stackrel{p}{\to}N\xrightarrow{\omega (p')^{-1}} M[1]$ is isomorphic to the distinguished triangle $\eta;$ 
proving that $M\stackrel{i}{\to} E \stackrel{p}{\to}N$ belongs to $\Ex\,(\Theta).$ Thus 
$\Ex_Q\,(\Theta)=\Ex\,(\Theta).$
\

(b) By (a), it follows that the derived functor $$R\Hom_\T(Q,-):\D^b(\F(\Theta))\to \D^b(\F({}_A\Delta))$$ is an equivalence of 
triangulated categories. Hence, (b) is a consequence of \ref{R1}.
\

(c) By \ref{ProAusl} (b) and \ref{exacto}, we get that $e_Q:\add\,(Q)\to \proj\,(A)$ is an exact equivalence of $R$-categories. Hence, we have that 
$K^{-,b}(\add\,(Q))\simeq K^{-,b}(\proj\,(A))$ as triangulated categories. Therefore (c) follows from (b), since 
$K^{-,b}(\proj\,(A))\simeq \D^b(A)$ as triangulated categories.
\end{dem}

\begin{teo}\label{R3} Let $(\Theta,\leq)$ be a $\Theta$-system, of size $t,$ in an artin triangulated $R$-category $\T.$ Then, the following statements hold.
\begin{itemize}
\item[(a)] The pair $(\F(\Theta),\Ex\,(\Theta))$ is an exact and Krull-Schmidt category.
\item[(b)] There exist a unique, up to isomorphism,  families $\Qb$ and $\Yb$ of objects 
in $\T,$ such that $(\Theta,\Qb,\leq)$ is a $\Theta$-projective system 
 and $(\Theta,\Yb,\leq)$ is a $\Theta$-injective system.
\item[(c)] For the $R$-algebras $A:=\End_\T(Q)^{op}$ and $B:=\End_\T(Y),$ the derived functors 
$$R\Hom_\T(Q,-):\D^b(\F(\Theta))\to \D^b(A)\;\text{ and }\;R\Hom_\T(-,Y):\D^b(\F(\Theta))\to \D^b(B)$$ are equivalences as triangulated categories.
\item[(d)] Both pairs $(A,\leq)$ and $(B,\leq^{op})$ are standardly stratified algebras, and 
moreover, the algebras $A$ and $B$ are derived equivalent.
\end{itemize}
\end{teo}
\begin{dem} Since the pair $(\Theta,\leq)$ is a $\Theta$-system, of size $t,$ in an 
artin triangulated $R$-category $\T,$ it follows from \ref{sisexproy} and its dual that (b) 
is true. Therefore, by \ref{R2} and its dual, we get (a) and (c). The fact that both pairs 
$(A,\leq)$ and $(B,\leq^{op})$ are standardly stratified algebras, can be obtained from 
\ref{TeoEquFilt} (c) and its dual. Finally, the fact that  $\D^b(A)\simeq\D^b(B)$ as triangulated 
categories (see (c)) say us that $A$ and $B$ are derived equivalent.
\end{dem}

\section{Examples}

\subsection{From stratifying systems in module categories} 
Let $\Lambda$ be an artin $R$-algebra and let $\D^b(\Lambda)$ be the bounded derived category of complexes in  $\modu\,(\Lambda).$ It 
is well-known that the canonical functor 
$\imath_0:\modu\,(\Lambda)\to\D^b(\Lambda),$ which sends $M\in\modu\,(\Lambda)$ to the 
stalk complex $M[0]$ concentrated in degree zero, is additive full and faithful. Hence, through 
the functor $\imath_0,$  the module category  $\modu\,(\Lambda)$ can be considered as a full 
additive subcategory of $\D^b(\Lambda).$ Furthermore 
$\Ext_\Lambda^k(X,Y)\simeq\Hom_{\D^b(\Lambda)}(X[0],Y[k])$ for any $k\in\Z$ and 
$X,Y\in\modu\,(\Lambda).$
\

In what follows, we recall from \cite{MMS1} the notion of stratifying systems; for a further development of such systems, see in \cite{ES,MMS1,MMS2,MSX}.

\begin{defi}\cite{MMS1} A stratifying system $(\Theta,\leq),$ of size $t$ in $\modu\,(\Lambda)$ consist of the 
following data.
\begin{enumerate}
\item[(SS1)] $\leq$ is a linear order on $[1,t].$
\item[(SS2)] $\Theta=\{\Theta(i)\}_{i=1}^{t}$ is a family of indecomposable objects in $\modu\,(\Lambda)$.
\item[(SS3)] $\Hom_{\Lambda}(\Theta(j),\Theta(i))=0$ for $j>i$.
\item[(SS4)] $\Ext_{\Lambda}^1(\Theta(j),\Theta(i))=0$ for $j\geq i$.
\end{enumerate}
\end{defi} 

By using the formula $\Ext_\Lambda^k(X,Y)\simeq\Hom_{\D^b(\Lambda)}(X[0],Y[k]),$ we get 
that any stratifying system $(\Theta,\leq),$ of size $t$ in $\modu\,(\Lambda),$ produces the 
$\Theta[0]$-system $(\Theta[0],\leq)$ of size $t$ in the triangulated category 
$\D^b(\Lambda).$

\subsection{Exceptional sequences} 

The notion of exceptional sequence originates from the study of vector bundles 
(see, for instance, \cite{Bo,GL}). Here, $\T$ denotes an artin triangulated $R$-category.

\begin{defi}\cite{Bo,Ru} An exceptional sequence of size $t,$ in the triangulated category $\T,$ is 
a sequence $\E=(\E_1,\E_2,\cdots,\E_t)$ of objects in $\T$ satisfying the following conditions.
\begin{enumerate}
\item[(ES1)] $\End_\T(\E_i)$ is a division ring, for each $i\in[1,t].$
\item[(ES2)] $\Hom_\T(\E_i,\E_i[k])=0$ $\forall\,i\in[1,t],$ $\forall\,k\in\Z-\{0\}.$
\item[(ES3)] $\Hom_\T(\E_j,\E_i[k])=0$ for $j>i$ and $\forall\,k\in\Z.$
\end{enumerate}
An exceptional sequence $\E=(\E_1,\E_2,\cdots,\E_t)$ is called strongly exceptional if the 
condition $\mathrm{(ES4)}$ holds, where 
$$\mathrm{(ES4)}\;\Hom_\T(\E_i,\E_j[k])=0\;\forall\,i,j\in[1,t],\;\forall\,k\in\Z-\{0\}.$$
\end{defi}

We recall that strongly exceptional sequences appear very often in algebraic geometry and provides a non-commutative model for the study of algebraic varieties (see \cite{Bo}).
\
 
Observe that any strongly exceptional sequence $\E=(\E_1,\E_2,\cdots,\E_t)$ of size $t,$ in the 
triangulated category $\T,$ is an example of a homological system in $\T.$ Namely, the pair 
$(\E,\leq),$ for $\leq$ the natural order on $[1,t],$ is an $\E$-system in $\T.$ So, as an application of the developed theory of homological systems we get the following result. 

\begin{teo}\label{TeoSES} Let $\E=(\E_1,\E_2,\cdots,\E_t)$ be a strongly exceptional sequence in 
an artin triangulated $R$-category $\T,$ and let $E:=\oplus_{i=1}^t\,\E_i.$ Then, the following statements hold.
\begin{itemize}
\item[(a)] The pair $(\F(\E),\Ex\,(\E))$ is an exact and Krull-Schmidt category.
\item[(b)] For the $R$-algebra $A:=\End_\T(E)^{op},$ the derived functor 
$$R\Hom_\T(E,-):\D^b(\F(\E))\to \D^b(A)$$ is an equivalence as triangulated categories.
\item[(c)] The pair $(A,\leq)$ is a quasi-hereditary algebra. 
\end{itemize}
\end{teo}
\begin{dem} By the condition $\mathrm{(ES4)},$ it follows that the triple $(\E,\E,\leq)$ is 
the $\E$-projective system associated to the $\E$-system $(\E,\leq).$ Thus, from 
\ref{R3} and the definition of strongly exceptional sequence the result follows.
\end{dem}
\begin{rk} Let $\T:=\D^b(\mathrm{Sh}(X))$ be the bounded derived category of coherent sheaves on a smooth manifold 
$X.$ 
\begin{itemize}
 \item[(1)] In \cite[Theorem 6.2]{Bo} it is proven that, for any strongly exceptional sequence $\E$ in $\T,$ such that $\T$ is 
 generated by $\E,$ the triangulated category $\T$ is equivalent to the bounded derived category $\D^b(A),$ where  
 $A:=\End_\T(E)^{op}.$
 \item[(2)] Observe that, in \ref{TeoSES}, it is not assumed that  $\T$ is generated by the strongly exceptional sequence $\E.$
\end{itemize}
\end{rk}

\vspace{.3cm}

Now, we consider a hereditary abelian $k$-category $\mathcal{H},$ for some field $k.$ By a result of 
Ringel (see \cite[Theorem 1]{R2}), it follows that $\Hom_\T(\Theta,\Theta[-1])=0$ for any set 
$\Theta$ of indecomposable objects in $\T:=\D^b(\mathcal{H}).$ Hence, we get that any 
exceptional sequence $\E=(\E_1,\E_2,\cdots,\E_t)$ of size $t,$ in the 
bounded derived category $\D^b(\mathcal{H}),$ is an example of a homological system in $\D^b(\mathcal{H}).$ Namely, the pair 
$(\E,\leq),$ for $\leq$ the natural order on $[1,t],$ is an $\E$-system in $\D^b(\mathcal{H}).$ So, as an application of the developed theory of homological systems we get the following result. 

\begin{teo}\label{TeoES} Let $\E=(\E_1,\E_2,\cdots,\E_t)$ be an exceptional sequence in the 
triangulated category $\T:=\D^b(\mathcal{H})$. Then, the following statements hold true.
\begin{itemize}
\item[(a)] The pair $(\F(\E),\Ex\,(\E))$ is an exact and Krull-Schmidt category.
\item[(b)] There exist a unique, up to isomorphism,  families $\Qb$ and $\Yb$ of objects 
in $\T,$ such that $(\E,\Qb,\leq)$ is a $\E$-projective system 
 and $(\E,\Yb,\leq)$ is a $\E$-injective system.
\item[(c)] For the $R$-algebras $A:=\End_\T(Q)^{op}$ and $B:=\End_\T(Y),$ the derived functors 
$$R\Hom_\T(Q,-):\D^b(\F(\E))\to \D^b(A)\;\text{ and }\;R\Hom_\T(-,Y):\D^b(\F(\E))\to \D^b(B)$$ are equivalences as triangulated categories.
\item[(d)] Both pairs $(A,\leq)$ and $(B,\leq^{op})$ are quasi-hereditary algebras, and 
moreover, the algebras $A$ and $B$ are derived equivalent. 
\end{itemize}
\end{teo}
\begin{dem} It follows from \ref{R3} and the definition 
of a exceptional sequence.
\end{dem}

\subsection{A $\Theta$-system which is not an exceptional sequence} 
In what follows, we give an example of a $\Theta$-system which is not a exceptional sequence and does not come from a stratifying system in a module category. To see that, we consider the hereditary path $k$-algebra 
$\Lambda:=k(1\fl 2\fl 3)$ and the triangulated category $\T:=\D^b(\Lambda).$ The Auslander-Reiten quiver of the bounded 
derived category $\D^b(\Lambda)$ can be seen in the Figure 1.
\

\begin{figure}[h]
\[ \scalebox{1}{
\begin{tikzpicture}[scale=0.8,>=stealth]
\node (P1) at (0,0) {$P_{1}[0]$};

\node (P2) at (1,1) {$P_{2}[0]$};

\node (P3) at (2,2) {$P_{3}[0]$};

\node (S2) at (2,0) {$S_{2}[0]$};

\node (I2) at (3,1) {$\bf{I_{2}[0]}$};

\node (I3) at (4,0) {$I_{3}[0]$};

\node (P11) at (4,2) {$P_{1}[1]$};

\node (P21) at (5,1) {$P_{2}[1]$};

\node (P31) at (6,0) {$P_{3}[1]$};

\node (S21) at (6,2) {$S_{2}[1]$};

\node (I21) at (7,1) {$I_{2}[1]$};

\node (I31) at (8,2) {$I_{3}[1]$};


\node (P12) at (8,0) {$P_{1}[2]$};

\node (P22) at (9,1) {$P_{2}[2]$};

\node (P32) at (10,2) {$P_{3}[2]$};

\node (S22) at (10,0) {$S_{2}[2]$};

\node (I22) at (11,1) {$\bf{I_{2}[2]}$};

\node (I32) at (12,0) {$I_{3}[2]$};







\node (I1-1) at (-1,1) {};

\node (P23) at (13,1) {};

\draw [->] (P1) -- (P2);

\draw [->] (P2) -- (P3);

\draw [->] (P2) -- (S2);

\draw [->] (P3) -- (I2);

\draw [->] (S2) -- (I2);

\draw [->] (I2) -- (I3);

\draw [->] (P11) -- (P21);

\draw [->] (P21) -- (P31);

\draw [->] (P21) -- (S21);

\draw [->] (P31) -- (I21);

\draw [->] (S21) -- (I21);

\draw [->] (I21) -- (I31);

\draw [->] (P12) -- (P22);

\draw [->] (P22) -- (P32);

\draw [->] (P22) -- (S22);

\draw [->] (P32) -- (I22);

\draw [->] (S22) -- (I22);

\draw [->] (I22) -- (I32);







\draw [loosely dotted, thick] (I1-1)--(P2);

\draw [loosely dotted, thick] (I22)--(P23);

\draw [->] (I2) -- (P11);

\draw [->] (I3) -- (P21);

\draw [->] (I21) -- (P12);

\draw [->] (I31) -- (P22);

\end{tikzpicture}
} \] \caption{\label{carcaj} The bounded derived category $\D^b(\Lambda)$.}
\end{figure}
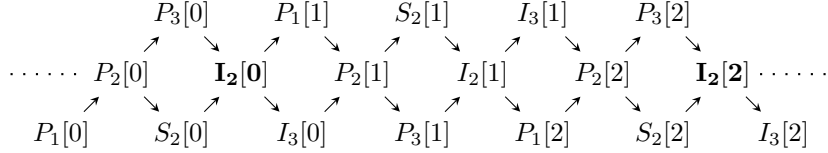

Consider the natural order $1\leq 2\leq 3$ and the set $\Theta:=\{\Theta(1),\Theta(2), \Theta(3)\},$ 
of indecomposable objects in $\T,$ where 
$\Theta(1):=I_{2}[0]$, $\Theta(2):=I_{2}[2]$ and $\Theta(3):=I_{2}[4]$. We assert that the pair 
$(\Theta,\leq)$ is a $\Theta$-system of size 3 in the triangulated category $\T.$
Indeed, by using Figure 1, it can be checked that $\Hom_\T(\Theta(j),\Theta(i))=0$ for $j >i$. The 
condition $\Hom_\T(\Theta,\Theta[-1])$ follows from the fact that $\modu\,(\Lambda)$ is an 
abelian hereditary $k$-category. Finally, using that 
$\Ext_{\Lambda}^{k}(X,Y)\simeq\Hom_T(X[0],Y[k]),$ it can be seen that $\Hom_\T(\Theta(j),\Theta(i)[1])=0$ 
for $j\geq i$. Therefore, the pair $(\Theta,\leq)$ is a $\Theta$-system in $\T.$
\

Observe that $\Hom_\T(\Theta(3),\Theta(2)[2])=\Hom_\T(\Theta(3),\Theta(3))\neq
0$, and so $\Theta$ is not an exceptional sequence. Furthermore the pair $(\Theta,\leq)$ does not 
come from a stratifying system in the module category $\modu\,(\Lambda).$
\vspace{.3cm}

{\bf Acknowledgments} Both authors would like to thank Mar\'ia Jose Souto Salorio for useful discussions and suggestions.

\footnotesize

\vskip3mm \noindent Octavio Mendoza Hern\'andez:\\ Instituto de Matem\'aticas, Universidad Nacional Aut\'onoma de M\'exico\\
Circuito Exterior, Ciudad Universitaria,
C.P. 04510, M\'exico, D.F. MEXICO.\\ {\tt omendoza@matem.unam.mx}

\vskip3mm \noindent Valente Santiago Vargas:\\ Instituto de Matem\'aticas, Universidad Nacional Aut\'onoma de M\'exico\\
Circuito Exterior, Ciudad Universitaria,
C.P. 04510, M\'exico, D.F. MEXICO.\\ {\tt valente@matem.unam.mx}

\end{document}